\documentclass{amsart} 
\usepackage{fullpage}
\usepackage[utf8]{inputenc}
\usepackage{dsfont}
\usepackage[all,2cell,cmtip]{xy}
\UseTwocells
\usepackage{amssymb}
\usepackage{amsmath,pdftexcmds,mathtools}
\usepackage[T1]{fontenc}

\let\oldtocsection=\tocsection

\let\oldtocsubsection=\tocsubsection

\let\oldtocsubsubsection=\tocsubsubsection

\renewcommand{\tocsection}[2]{\hspace{0em}\oldtocsection{#1}{#2}}
\renewcommand{\tocsubsection}[2]{\hspace{1em}\oldtocsubsection{#1}{#2}}
\renewcommand{\tocsubsubsection}[2]{\hspace{2em}\oldtocsubsubsection{#1}{#2}}

\usepackage{epigraph}

\DeclareFontFamily{U}{cbgreek}{}
\DeclareFontShape{U}{cbgreek}{m}{n}{
        <-6>    grmn0500
        <6-7>   grmn0600
        <7-8>   grmn0700
        <8-9>   grmn0800
        <9-10>  grmn0900
        <10-12> grmn1000
        <12-17> grmn1200
        <17->   grmn1728
      }{}
\DeclareFontShape{U}{cbgreek}{bx}{n}{
        <-6>    grxn0500
        <6-7>   grxn0600
        <7-8>   grxn0700
        <8-9>   grxn0800
        <9-10>  grxn0900
        <10-12> grxn1000
        <12-17> grxn1200
        <17->   grxn1728
      }{}

\makeatletter
\newcommand{\normalorbold}{%
  \ifnum\pdf@strcmp{\math@version}{bold}=\z@ bx\else m\fi
}
\makeatother

\usepackage{graphicx}
\usepackage{tikz}
\usepackage{array}

\usepackage{todonotes}

\allowdisplaybreaks[3]

\usepackage{epigraph}

\usetikzlibrary{decorations.pathmorphing}

\usepackage[outline]{contour}
\contourlength{2pt}

\usepackage{bigdelim}
\usepackage{enumerate}
\usepackage{mathrsfs}
\usepackage{color}
\usepackage{amsthm}
\usepackage{float}

\usepackage{accents}

\usepackage{stmaryrd}

\usepackage{physics}

\usepackage{tensor}

\theoremstyle{plain}
\newtheorem{theorem}{Theorem}[section]
\newtheorem{proposition}[theorem]{Proposition}

\newtheorem{lemma}[theorem]{Lemma}
\newtheorem*{theorem*}{Theorem}

\theoremstyle{definition}
\newtheorem{definition}[theorem]{Definition}

\usepackage{txfonts}

\usepackage{multirow}
\usepackage{multicol}
\title{Operator $K$-theory algebra spectra of $C^*$-algebras}
\author{R. Vasconcellos V. \and L. C. P. A. M. Müssnich \and N. J. B. Aza}

\begin{document}

\begin{abstract}
We construct commutative algebra spectra that represent the operator $K$-theory of $C^*$-algebras, which are algebras over the commutative ring spectra that represent topological $K$-theory. The spectral multiplicative structure introduces a new graded commutative ring structure on the $K$-groups, generalizing the well-known graded ring structure of commutative $C^*$-algebras. This last structure reflects the multiplicative structure of topological $K$-theory via Gelfand duality, Swan's theorem and the fiber tensor product.

We introduce $\mathscr L$-permutative categories, a generalization of bipermutative categories, which are permutative categories equipped with a multiplicative structures induced by coherent actions of the linear isometries operad. The main class of examples of interest are categories whose objects are projection matrices of the unitization of the stabilizations $\widetilde{\mathfrak{KA}}$ of a $C^*$-algebras $\mathfrak A$, and morphisms partial isometries witnessing the Murray-von Neumann relation.

We then construct $E_\infty$-ring spaces out of them by adapting the usual method applied to bipermutative categories. The delooping functor of the recognition principle, the homotopical augmentation ideal and localization at the Bott element then give us our algebra spectra.
\end{abstract}

\maketitle

\begin{center}
\begin{multicols}{2}
\footnotesize \raggedright
    A physicist is flying in a balloon
    
    caught in a great storm
    
    thrown here and there
    
    blinded by rain and tears
    
    \vspace{0.25cm}
    
    As the storm subsides
    
    the dawn is visible on the horizon
    
    the physicist overflies a desert
    
    a single person walks beneath
    
    \vspace{0.25cm}
    
    ``Where am I?''
    
    shouts the physicist
    
    \vspace{0.25cm}
    
    The lone walker looks up
    
    and waves
    
    \vspace{0.25cm}
    
    They enter a state of deep thought
    
    \vspace{0.25cm}
    
    As the sun sets
    
    slowly waking
    
    the physicist hears an enthusiastic voice:
    
    ``You are in a balloon!!!''
    
    \vspace{0.25cm}
    
    ...
    
    \vspace{0.25cm}
    
    ``Are you a mathematician?''
    
    \vspace{0.25cm}
    
    ``Yes
    
    how do you know?''
    
    \vspace{0.25cm}
    
    ``I asked a simple question
    
    you took a long time to give me an answer
    
    which was completely correct
    
    and it didn't help me at all''
\end{multicols}
\end{center}
\hspace{1.3cm}Our hope for this article is that a detailed description of the balloon might help.

\tableofcontents

\section{Introduction}
Additive topological $K$-theory extends to a homology theory $K_-$ for $C^*$-algebras, called operator $K$-theory \cite{BlackadarOpKTh,schroder1993k}. For a unital $C^*$-algebra $\mathfrak A$ consider the groupoid $\texttt{pr}_{\mathfrak A}$ of projection matrices, i.e. self-adjoint idempotent matrices, and partial isometries witnessing the Murray-von Neumann relation between them. The direct sum of matrices gives us a permutative category structure on $\texttt{pr}_{\mathfrak A}$, that is a strictly associative and unital symmetric monoidal structure. The categories $\texttt{pr}_{\mathfrak A}$ and $\texttt{prmod}_{\mathfrak A}$ of projective finitely generated right $\mathfrak A$-modules are related by the fact that the images of the projections matrices are projective right $\mathfrak A$-modules, with the direct sum of projection matrices mapping to the direct sum of the respective modules.

If $\mathfrak A$ is commutative then the Kronecker product induces a bipermutative structure on $\texttt{pr}_{\mathfrak A}$, with the multiplicative structure mapping to the tensor product of  modules. The Kronecker product is ill-behaved when $\mathfrak A$ is noncommutative, which is related to the lack of a useful notion of tensor product of modules over noncommutative rings. This means that the multiplicative structure of topological $K$-theory doesn't extend naively to the noncommutative setting.

Our main goal is to establish a multiplicative structure for operator $K$-theory by constructing for each $C^*$-algebra $\mathfrak A$ a ring spectrum whose stable homotopy groups are naturally isomorphic to the $K$-groups of $\mathfrak A$. We will consider both complex and real operator $K$-theory. For $\mathds F$ equal to either $\mathds R$ or $\mathds C$ and $\mathfrak A$ a $C^*$-algebra over $\mathds F$ the functorialy constructed commutative ring spectrum $\widetilde K\mathfrak A$ will be a $K\mathcal U_{\mathds F}$-algebra spectrum, where $K\mathcal U_{\mathds F}$ is the commutative ring spectrum representing the topological $K$-theory over $\mathds F$. Our main theorem can be stated as follows:
\begin{theorem*}
    There is a functor 
    $$\widetilde K\in\texttt{Cat}(
    \texttt{C$^*$Alg}_{\mathds F},
    \texttt{AlgSp}_{K\mathcal U_{\mathds F},\text{nu}})
    $$
    such that the induced homology theory $\pi^S_-\widetilde K$ is naturally isomorphic to $K_-$.
    
    The spectral multiplicative structure of $\widetilde K\mathfrak A$ induces a natural graded-commutative ring structure on $K_-\mathfrak A$.
\end{theorem*}

Constructions of spectra that represent certain versions of operator $K$-theory can be found in \cite{BJS,DEKM}. Multiplicative $K$-theoretical structures were introduced for $C^*$-algebras associated to Smale spaces in \cite{Put96}, and operator $K$-theory commutative ring spectra were constructed for the subcategory of separable, unital and strongly self-absorbing $C^*$-algebras \cite{DP15}. Our construction can be applied to any $C^*$-algebra, and is functorial on $*$-homomorphisms. To the best of our knowledge there is no other construction of multiplicative operator $K$-theory with this level of generality.

The spectrum $K \mathcal U_{\mathds F}$ can be constructed from the delooping of the $E_\infty$-ring space associated to the  bipermutative category $\mathcal U_{\mathds F}$ of unitary matrices by localizing at the appropriate Bott element \cite{MayMultLoopSpcTh,MaEinftyRingSpcsSpectra,Ma09b}. Our method is an adaptation of this construction, where we substitute bipermutative categories with categories equipped with operad theoretical structures we introduce. This construction goes through the category of $E_\infty$-ring spaces, which suggests it might have orientation theoretical applications in the theory of noncommutative bundles and fibrations, similar to the ones reviewed in \cite{Ma09c} in the commutative setting.

The first step of our construction is to consider for a $C^*$-algebra $\mathfrak A$ the unitization of its stabilization $\widetilde{\mathfrak{KA}}$. Taking the unitization is a standard step in the construction of operator $K$-theory which guarantees half-exactness, while taking the stabilization is well known to preserve the additive structure and is essential for our construction of the multiplicative structure.
If $\mathfrak A$ is commutative the Kronecker product preserves projections and partial isometries in $\texttt{pr}_{\widetilde{\mathfrak{KA}}}$, which isn't the case if $\mathfrak A$ is noncommutative. Our strategy here is to use the linear isometries operad $\mathscr L$ to parameterize a family of Kronecker-like products on $\texttt{pr}_{\widetilde{\mathfrak{KA}}}$ that do preserve projections and partial isometries, and are also compatible with the direct sum. Thus our construction provides $\mathscr L$-parameterized families of representatives for products, even if no canonical representatives exist.

We formalize the algebraic structure of $\texttt{pr}_{\widetilde{\mathfrak{KA}}}$ by generalizing the definition of bipermutative categories through categorial actions of topological operad pairs. Our definition of $\mathscr G$-categories is different than that of $\mathcal A_\infty$-categories, such as Fukaya categories, which are weak-categories where compositions are parameterized by an $\mathcal A_\infty$-operad. In our definition the operad action is a structure separate, though compatible, with the composition.
In \cite{MayMultLoopSpcTh,Ma09b} it is shown how to construct $E_\infty$-ring spaces out of bipermutative categories. An $E_\infty$-ring space is like a topological semiring, except its operations satisfy the axioms only up to coherent homotopies. We construct $E_\infty$-ring spaces $\Lambda\lVert\mathfrak A\rVert$ out of the $ \texttt{pr}_{\widetilde{\mathfrak{KA}}}$ by essentially the same method. The only departure from the original construction is that we begin by defining lax functors with domain the category of ring operators $\hat{\mathscr L}\wr\mathscr F$, instead of $\mathscr F\wr\mathscr F$. From this point on the classical recipe takes over. Applying the multiplicative recognition theorem of \cite{MaEinftyRingSpcsSpectra,MayMultLoopSpcTh,Ma09a} and homotopy fibers we obtain a connective commutative $k\mathcal U_{\mathds F}$-algebra spectrum $\widetilde k\mathfrak A$ with the correct non-negative stable homotopy groups. We further localize on Bott elements in order to construct a $K\mathcal U_{\mathds F}$-algebra spectrum $\widetilde K \mathfrak A$ with the right negative stable homotopy groups.

The multiplicative structure of topological $K$-theory has many important consequences (see for instance \cite[II.2.6]{atiyah2018k}). In general having an algebraic structure described by spectra is useful since a lot of classical algebraic theories and results generalize to this setting. For instance there is a good notion of Galois Theory for ring spectra \cite{rognes2008galois}. We can also apply some modern machinery such as topological Hochschild homology of algebra spectra \cite{EKMM}. 
Also, though we do not cover any equivariant $K$-theory considerations, the construction presented 
here should admit a similar extension. For this topic, we recommend \cite{gaunce2006equivariant}.

Apart from the interest that the above-mentioned multiplicative
structure might elicit within the context of K-theory for $C^{*}$--algebras, 
we expect that it should be relevant in the context of mathematical-physics. 
As an example, as it is widely known, K-theory for $C^{*}$-algebras
have a myriad of applications for instance in the context of condensed matter physics.
For instance, in the analysis of topological insulators \cite{cornfeld2021tenfold, prodan2016bulk}, 
the ring spectral structure of topological $K$-theory in particular plays an important role in the 
classification of states in crystalline topological insulators and superconductors.
To go even further, in quantum many--body systems, it is of plethora interest to classify 
\emph{free fermion} models by considering the underlying physical symmetries. In \cite{cornfeld2021tenfold}, 
and references therein, it is stressed the importance of real operator $K$-theory -- as well as the use 
of equivariant ring spectra -- to understand this physical situation. In a broad sense, for any spatial 
dimension $d\in\mathds{N}$, the mathematical and physical interpretation is summarized in the \emph{periodic} 
table of topological and superconductors, TISCs. In particular, it is assumed that the systems are in \emph{zero} 
temperature and that the unique relevant fermion-fermion interaction is provided by \emph{Pauli's principle}. 
Nonetheless, to the best of our knowledge, there is not a rigorous result for the interparticle case yet, where 
\emph{non-null} fermionic interactions -- or, at least, \emph{small enough} non-null interactions --
are taken into account

Concerning the latter, there is an increasing interest, within the mathematical-physics community,
to combine the already successful use of $K$--theory (used in the free fermion case, as stressed before) 
with $C^{*}$--dynamical systems, which is the accepted mathematical framework when one considers the 
interparticle setting \cite{bru_pedra16}. In this regard, for the \emph{universal} unital 
$\mathrm{CAR}$ $C^{*}$--algebra $\mathfrak{F}$, its associated set of \emph{states} 
$\mathfrak{E}$\footnote{\label{foot:state} The linear functional $\omega\in\mathfrak{F}^{*}$ is named 
a state if it is positive and normalized: $\omega(A^{*}A)\geq 0$, for any $A\in\mathfrak{F}$, and 
$\omega(\tilde{1})=1$. In particular, $\overline{\omega\left(A\right)}=\omega\left(A^{*}\right)$, for any
$A\in\mathfrak{F}$. The set of all states on $\mathfrak{F}$ is denoted by $\mathfrak{E}\subset\mathfrak{F}^{*}$.}, 
and non--negatives temperatures $T\in\mathds{R}_{0}^{+}$, we have already studied some mathematical aspects of
interacting fermions \cite{LD1}. In particular, for \emph{gapped} fermions embedded in zero temperature systems, 
as is the case in the TISCs framework, we are able to obtain suitable estimates for well--defined infinite systems.
Moreover, in \cite{ARS22}, in order to distinguish between two gapped free fermion systems at zero temperature, a
$\mathds{Z}_{2}$--topological index was studied. Particularly, connectedness properties of the set of \emph{ground states} 
associated with our systems were considered. In \cite{AMR22}, the mathematical 
framework to generalize the main results given in \cite{ARS22} for the weakly interacting fermion setting
is stated. See the canonical references \cite{Araki-Moriya,araki1985ground,BratteliRobinson} for concrete 
definitions of ground states. For recent results of classification of interacting fermion systems and other of their results see
\cite{nach21quasi,ogata2021class2,ogata21cl}, and references therein. 

Therefore, given all this discussion, we strongly believe that a product for K-groups
will bring in a new tool for analyzing physical systems from a mathematical view-point.

\subsection{Structure of the article}

In the second section the basics of operator $K$-theory is presented, including the matrix operations we will need and their properties. This includes the construction of the Murray-von Neumann groupoids $\texttt{pr}_{\mathfrak A}$ of $C^*$-algebras and their additive permutative structure given by direct sums of matrices.

In the third section we introduce actions of operads on categories, and show that $\mathscr L$ acts on the Murray-von Neumann groupoids $\texttt{pr}_{\widetilde{\mathfrak{KA}}}$ of unitizations of stabilizations of $C^*$-algebras. We further define actions of operad pairs on categories, which encode coherent distributivity structures. We will be particularly interested in $\mathscr L$-permutative categories, which are $(\mathscr C,\mathscr L)$-categories, where $\mathscr C$ is the terminal operad. This is a generalization of bipermutative categories. This section is the conceptual heart of the article and contains the main new ideas we are introducing. 

The following sections are then adaptations of the construction of ring spectra that represent topological $K$-theory from the bipermutative categories $\mathcal U_{\mathds F}$. We leave most of the details and proofs of statements to the references, but we make an effort to explicitly lay out formulas for the algebraic structures involved.

The fourth section presents the construction of special $\hat{\mathscr L}\wr \mathscr F$-spaces out of $\mathscr L$-permutative categories by Street's rectifications of lax functors \cite{street1972two}.

In the fifth section we show how to apply the two-sided bar construction to construct $E_\infty$-ring spaces out of the $\hat{\mathscr L}\wr \mathscr F$-spaces of the previous section. In this and the previous section we follow closely the construction of $E_\infty$-ring spaces out bipermutative categories in \cite{MayMultLoopSpcTh,Ma09b}.

In the sixth section we apply the delooping functor of the multiplicative recognition theorem, homotopy fibers and localization at Bott elements to obtain $K\mathcal U_{\mathds F}$-algebra spectra that represent operator $K$-theory, similar to how $K\mathcal U_{\mathds F}$ is constructed in \cite{MaEinftyRingSpcsSpectra,Ma09a}.

\subsection{Notation and terminology}

Given a class $A$ and a family of classes $\langle B^a\rangle$ indexed by $A$ the \textit{dependent sum} $\Sigma_AB^a$ is the class of pairs $(a,b)$ with $a\in A$ and $b\in B^a$; the \textit{dependent product} $\Pi_AB^a$ is the class of sequences $\langle b^a\rangle$ indexed on $A$ with $b^a\in B^a$ for each $a\in A$, or equivalently it is the class of sections of the natural surjection $\Sigma_AB^a\rightarrow A$. If the $B^a$ equal a fixed $B$ then $\Sigma_A B=A\times B$ and $\Pi_A B=B^A$. For $A$ a set (or space) equipped with an equivalence relation $\sim$ we will denote the equivalence classes of $a\in A$ using square brackets $[a]\in A/_\sim$. We will express elements of mapping sets $B^A$ as $a\mapsto \Phi$ for some expression $\Phi$ which may use the variable $a$.

Let $\mathbb F$ be the category with set of objects the isomorphism classes of linearly ordered finite sets and morphisms the sets of (not necessarily monotonous) functions between them. All isomorphism classes in $\mathbb F$ can be represented by $\underline p=\{1<2<\cdots<p\}$ for some $p\in\mathds N$. We include in $\mathbb F$ the empty set $\underline 0:=\emptyset$. Let $\mathbb F^{\text{inj}}$ be the category with $\text{Ob }\mathbb F^{\text{inj}}:=\mathbb F$ and not necessarily monotonous injective functions as morphisms, and $\mathbb F^{\text{bij}}\subset \mathbb F^{\text{inj}}$ the subcategory with the same objects and with bijections as morphisms. For $(\underline a,\langle \underline p^x\rangle)\in \Sigma_{\mathbb F}\mathbb F^{\underline a}$ the dependent sum can be equipped with the lexicographical order, so that $\Sigma_{\underline a}\underline p^x\in\mathbb F$. Similarly for $(\underline m,\langle \underline p^i\rangle)\in \Sigma_{\mathbb F}\mathbb F^{\underline m}$ we also have $\Pi_{\underline m}\underline p^i\in \mathbb F$, equipped with the lexicographical order.

For all $(\underline a^0, \underline a^1,\sigma,\langle \underline p^{x_0}\rangle)\in\Sigma_{\mathbb F^2} \mathbb F(\underline a^0,\underline a^1)\times\mathbb F^{\underline a^0}$ we define
\begin{equation*}
    \sigma(\underline p^{x_0})\in\mathbb F(\Sigma_{\underline a^0}\underline p^{x_0},\Sigma_{\underline a^1}\underline p^{\sigma^{-1} x_1}), \qquad
    \sigma(\underline p^{x_0})(x'_0,\mu):=(\sigma x'_0,\mu),
\end{equation*}
where if $x_1\not\in \text{Im}\sigma$ then $\underline p^{\sigma^{-1} x_1}:=\underline 0$.

For all $(\underline a,\langle (\underline p^{x0},\underline p^{x1},\sigma^x)\rangle )\in \Sigma_{\mathbb F}\Pi_{\underline a}\Sigma_{\mathbb F^2}\mathbb F(\underline p^{x0},\underline p^{x1})$ we define
\begin{equation*}
    \Sigma_{\underline a}\sigma^x\in \mathbb F(\Sigma_{\underline a}\underline p^{x0},\Sigma_{\underline a}\underline p^{x1}),
    \qquad
    \Sigma_{\underline a}\tau ^x(x',\mu_0)=(x',\sigma^{x'}\mu_0).
\end{equation*}

For all $(\underline m^0, \underline m^1,\sigma,\langle \underline p^{i_0}\rangle)\in\Sigma_{\mathbb F^2} \mathbb F(\underline m^0,\underline m^1)\times\mathbb F^{\underline m^0}$ we define
\begin{equation*}
    \sigma\langle\underline p^{i_0}\rangle\in\mathbb F(\Pi_{\underline m^0}\underline p^{i_0},\Pi_{\underline m^1}\underline p^{\sigma^{-1} i_1}), \qquad
    \sigma\langle\underline p^{i_0}\rangle\langle \mu^{i'_0}\rangle:=\langle \mu^{\sigma^{-1} i_1}\rangle,
\end{equation*}
where if $i_1\not\in \text{Im}\sigma$ then $\underline p^{\sigma^{-1} i_1}:=\underline 1$.

For all $(\underline m,\langle (\underline p^{i0},\underline p^{i1},\sigma ^i)\rangle )\in \Sigma_{\mathbb F}\Pi_{\underline m}\Sigma_{\mathbb F^2}\mathbb F(\underline p^{i0},\underline p^{i1})$ we define
\begin{equation*}
    \Pi_{\underline m}\sigma^i\in \mathbb F(\Pi_{\underline m}\underline p^{i0},\Pi_{\underline m}\underline p^{i1}),
    \qquad
    \Pi_{\underline m}\sigma^i\langle \mu^{i'}_0\rangle=\langle \sigma^{i'}\mu^{i'}_0\rangle.
\end{equation*}

For all $(\underline m,\langle \underline a^i\rangle,\langle \underline p^{ix}\rangle)\in \Sigma_{\Sigma_{\mathbb F}\mathbb F^{\underline m}}\mathbb F^{\Sigma_{\underline m}\underline a^i}$ we define
$$
    \delta_{\underline m,\langle \underline a^i\rangle,\langle \underline p^{ix}\rangle}\in\mathbb F^{\text{bij}}(\Pi_{\underline m}\Sigma_{\underline a^i}\underline p^{ix},\Sigma_{\Pi_{\underline m}\underline a^i}\Pi_{\underline m}\underline p^{ix^i}),
    \qquad
    \delta_{\underline m,\langle \underline a^i\rangle,\langle \underline p^{ix}\rangle}\langle (x^{i'},\mu^{i'})\rangle:=(\langle x^{i'}\rangle,\langle \mu^{i'}\rangle),
$$
which in general are not monotonous.

Let $\mathscr F$ be the category of isomorphism classes of linearly ordered finite based sets and not-necessarily monotonous based functions. All isomorphism classes in $\mathscr F$ can be represented by $\underline m_*=(\{0<1<\cdots<m\},0)$ for some $m\in\mathds N$. Dependent sums induce dependent wedge sums $\vee_{\underline m}\underline p^i_*:=(\Sigma_{\underline m}\underline p^i)_*$ and the dependent product induces dependent smash products $\wedge_{\underline m}\underline p^i_*:=(\Pi_{\underline m}\underline p^i)_*$.

We denote by $\text{Gr}:\texttt{CMon}\rightarrow \texttt{AbGrp}$ \textit{Grothendieck's enveloping group of commutative monoids construction}. We also denote by $\text{Gr}:\texttt{SemiRing}\rightarrow \texttt{Ring}$ the \textit{enveloping ring of semirings construction}.

We denote by $\texttt{Top}$ category of compactly generated and weakly Hausdorff topological spaces, by $\texttt{Top}_{\text{LCH}}\subset \texttt{Top}$ the subcategory of locally compact Hausdorff spaces and proper maps, and by $\texttt{Top}_{\text{CH}}\subset \texttt{Top}$ the full subcategory of compact Hausdorff spaces. For $X\in\texttt{Top}_{\text{LCH}}$ we denote its one point compactification as $X_+\in \texttt{Top}_{\text{CH}}$. We denote by $I:=[0,1]\in\texttt{Top}_{\text{CH}}$ the real interval. 

We denote by $*\texttt{Ring}$ the category of \textit{$*$-rings}, i.e. rings equipped with an antiautomorphic involution $-^*$, and $*$-homomorphisms. Rings will not be assumed unital unless stated otherwise, and ring homomorphisms are not assumed to preserve units. A \textit{$*$-algebra} $\mathfrak A$ over a commutative $*$-ring $\mathds F$ is an associative $\mathds F$-algebra satisfying $(\alpha A)^*=\alpha ^*A^*$ for all $\alpha\in \mathds F$ and $A\in \mathfrak A$. The full subcategory of unital rings is denoted $*\texttt{Ring}_1$.

To any $*$-algebra $\mathfrak A$ over a unital $*$-ring $\mathds F$ we associate its \textit{unitization} $\widetilde{\mathfrak A}\in *\texttt{Ring}_1$ with underlying abelian group $\mathds F\oplus \mathfrak A$ and multiplication and involution expressed, using the notation $\alpha\widetilde 1+A:=(\alpha,A)\in \widetilde{\mathfrak A}$, as
$$(\alpha^1\tilde 1+A)(\alpha^2\tilde 1+A^2)
:=\alpha^1\alpha^2\tilde 1+\alpha^2 A^1+\alpha^1 A^2 +A^1A^2,
\qquad
(\alpha\tilde 1+A)^*:=(\alpha^*\tilde 1+A^*).
$$

A \textit{complex $C^*$-algebra} is a complex Banach $\ast$-algebra satisfying the $C^*$-identity $\lVert A^*A\rVert=\lVert A\rVert^2$. A \textit{real $C^*$-algebra} is a real Banach $\ast$-algebra satisfying the $C^*$-identity and the further condition that $\widetilde 1+A^*A\in \widetilde{\mathfrak A}$ is invertible for all $A\in \mathfrak A$. For $\mathds F$ equal to either $\mathds C$ or $\mathds R$ the category of $C^*$-algebras over $\mathds F$ and $*$-homomorphisms is denoted $\texttt C^*\texttt{Alg}_{\mathds F}$. The full subcategory of unital $C^*$-algebras is denoted as $\texttt C^*\texttt{Alg}_{\mathds F,1}$.

A \textit{topological category} is a category whose collections of objects and morphisms form topological spaces, and all structural maps are continuous. This is a stronger condition than being topologically enriched.

For $\mathcal C$ a category and $X\in \mathcal C$ we denote by $\mathcal C_{/X}$ the category of morphisms in $\mathcal C$ with codomain $X$ and commutative triangles as morphisms.

A \textit{permutative category} 
$\mathcal C=(\mathcal C,\circledast,\boldsymbol e,\tau^\circledast_{A,B})$ is a symmetric monoidal topological category in which associativity and identity holds strictly. Permutative categories and symmetric monoidal functors form a category $\texttt{PermCat}$. A \textit{bipermutative category} $\mathcal C$ is a topological category equipped with two permutative structures $(\mathcal C;\oplus,\boldsymbol 0,\tau^\oplus_{A,B})$ and $(\mathcal C;\otimes,\boldsymbol 1,\tau^\otimes_{A,B})$ plus a left distributivity natural transformation $\delta_{A,B,C}:A\otimes(B\oplus C)\Rightarrow (A\otimes B)\oplus(A\otimes C)$ such that $\otimes$ distributes strictly over $\oplus$ on the right, $\boldsymbol 0$ is a strict two-sided absorbing object for $\otimes$ and such that some diagrammatically determined natural coherence laws hold \cite{Ma09b}. Bipermutative categories and symmetric bimonoidal functors form a category $\texttt{BipCat}$.

Let $\texttt{Inn}$ denote the topologically enriched category of finite or countably infinite dimensional real inner product spaces and linear maps, with all spaces topologized as the colimit of the finite dimensional sub-spaces. Let $\mathscr I$ be the subcategory with the same objects and with morphisms the linear isometries. Both $\texttt{Inn}$ and $\mathscr I$ are monoidal under direct sums. For $\mathds U\in\texttt{Inn}$ we denote by $\mathscr A_{\mathds U}$ the set of finite dimensional subspaces of $\mathds U$, partially ordered by inclusion, and for $U\in\mathscr A_{\mathds U}$ we define $\mathscr A_U:=\{V\in\mathscr A_{\mathds U}\mid U\leq V\}$. For $\mathds U=\mathds R^\infty$ we simply write $\mathscr A:=\mathscr A_{\mathds R^\infty}$. For $\langle f_a\rangle\in\mathscr I(\oplus_A\mathds U^x,\mathds V)$ and $\langle \vec U^x\rangle \in \oplus_A\mathds U^x$ we use the Einstein summation convention $f_a\vec U^x:=\sum_Af_a\vec U^x$. For $\mathds U\in\mathscr I$ and $U< \mathds U$ we use the notation $U^\bot:=\{\vec v\in\mathds U\mid \forall \vec u\in U:\vec v\cdot\vec u=0\}$ for the orthogonal complement. 

For any $U< \mathds U$ we define $\text{pr}_U\in \texttt{Inn}(\mathds U,\mathds U)$ as the orthogonal projection onto $U$, and for $f\in\texttt{Inn}(\mathds U,\mathds V)$ define $f\!\!\restriction_U\in \texttt{Inn}( U,\mathds V)$ as the restriction of $f$ on $U$. For $f\in\mathscr I(\mathds U,\mathds V)$ its adjoint is $f^*:=f^{-1}\text{pr}_{f\mathds U}\in \texttt{Inn}(\mathds V,\mathds U)$. For all $U\in \mathscr A_{\mathds U}$ let $\mathds S^U$ be the one point compactification of $U$ obtained by adding a point $\infty$ at infinity and for $(U,V)\in \Sigma_{\mathscr A}\mathscr A_U$ let $V-U:=V\cap U^\bot$.

Let $\Delta$ be the simplicial category, whose set of objects are the linearly ordered sets $\underline s_*=\{0<1<\cdots<s\}$ and the morphism the monotonous functions. The cosimplicial space of partitions of the interval $\text{Part}\in\texttt{Top}^{\Delta}$, with each $\text{Part} \underline s_*$ topologized as a subspace of monotonous functions from $\underline{s}=\{1<\cdots<s\}$ to $I$ (see \cite[Section 1.2]{Vi21} for a description of the face and degeneracy maps). For each $\langle\langle t^{id}\rangle\rangle\in\Pi_{\underline m}\text{Part}\underline{s}^i_*$ the order of the points $t^{id}$ in $I$ induces an order on $\vee_{\underline m}\underline s^i_*$, and so an element
$$
    \textstyle\lhd_{\underline m} \langle t^{id}\rangle\in \text{Part}\vee_{\underline m}\underline{s}^i_*.
$$
For each $\langle\langle  t^{id}\rangle\rangle\in \Pi_{\underline m}\text{Part} \underline s ^{i}_*$ and $i'\in \underline m$ we can define
$$
    \textstyle \delta^{i'}\in \Delta(\vee_{\underline m}\underline{s}^{i}_*,\underline s^{i'}_*
    ),
    \qquad
    \delta^{i'}(i,d):=\begin{cases}
        0&t^{id}<t^{i'1}\\
        \max_{t^{i'd'} \leq t^{id}}i',
        &t^{id}\geq t^{i'1}
    \end{cases}
$$
such that $\delta^{i'}\cdot\lhd_{\underline m} \langle t^{id}\rangle = \langle t^{i' d}\rangle\in \text{Part}\underline s^{i'}_*$.

For any simplicial space $X\in\texttt{Top}^{\Delta^{\text{op}}}$ its geometric realization $|X|$ is defined via the coend construction \cite{Lo15} as
$$
    |X|:=\int^\Delta X\underline s_*\times \text{Part }\underline s_*.
$$

The reason we consider the geometric realization via the partitions cosimplicial space instead of the usual homeomorphic cosimplicial space of topological simplexes is that this choice simplifies the algorithm in \cite[Theorem 11.5]{may2006geometry}.

The nerve of a category is $\mathcal{NC}$ with $\mathcal{NC}\underline s_*:=\texttt{Cat}(\underline s_*,\mathcal C)$, where here we think of the order structure of $\underline s_*$ as a category structure. The geometric realization of a category is $\lvert\mathcal{C}\rvert:=\lvert \mathcal{NC}\rvert$.

The main tool for constructing spaces out of monadic inputs is the bar construction, which are the geometric realizations of natural resolutions of the monad algebras \cite[Section 9]{may2006geometry}. Let $\mathcal T$ and $\mathcal A$ be topological categories. The category $B(\mathcal T,\mathcal A)$ is defined as follows: the objects of $B(\mathcal T,\mathcal A)$ are triples $(F,C,X)$, with $(C,\mu,\eta)$ a monad in $\mathcal T$, $(F,\lambda)$ a $C$-functor in $\mathcal A$ and $(X,\xi)$ a $C$-algebra, and morphisms spaces are composed of triples $(\alpha,\phi,f)$ with $\phi$ a monad morphism, $f$ a $C$-morphism and $\alpha$ a $C$-functors morphism. The \textit{two-sided bar construction} is the functor
\begin{gather*}
    B_-:B(\mathcal T,\mathcal A)\rightarrow \mathcal A^{\Delta^\text{op}},
    \qquad
    B_{\underline s_*}(F,C,X):=FC^s X,\\
    s_d=FC^{i}\eta_{C^{s-d+1}},
    \quad \partial_i=
    \begin{cases}
        \lambda_{C^{s-1}},&d=0\\
        F C^{d-1}\mu_{C^{s-d}},&0<d<s\\
        F C^{s-1}\xi,&d=s
    \end{cases}
\end{gather*}

We will use the notation $B(F,C,X):=\lvert B_-(F,C,X)\rvert\in\texttt{Top}$ for the geometric realization of the bar construction.

\subsection{Acknowledgments}

R. Vasconcellos V. was financed by the grant 2020/06159-5, São Paulo Research Foundation (FAPESP).

\section{Operator {$K$}-theory}

We briefly review the basic $K$-theory of $C^*$-algebras we will need. Standard references to the theory are \cite{BlackadarOpKTh,schroder1993k}.

We begin by establishing notations and basic results for matrix operations. We will extensively use dependent sums and products of linearly ordered sets to index matrix entries, since this will help streamline our main proofs.

We then give a brief review of the basic definitions and results of operator $K$-theory that we will need.

\subsection{Matrix operations}

Let $R\in *\texttt{Ring}$ and $\underline p^0,\underline p^1\in \mathbb F$. We denote by $ M_{\underline p^0,\underline p^1}R$ the set of $p^0\times p^1$ matrices $U=\mqty[U_{\mu\nu}]$. We set the convention that $ M_{\underline p^0,\underline p^1}R:=\left\{0_{\underline 0}\right\}$ whenever $p^0=0$ or $p^1=0$, where $0_{\underline 0}$ is a formal empty matrix. We also denote the union of all matrices over $R$ as $ MR:=\coprod_{\mathbb F^2} M_{\underline p^0,\underline p^1}R$.

We will extensively use Einstein's summation convention for matrix products, i.e.
$$ (U^1,U^2)\in M_{\underline p^0,\underline p^{1}}R\times M_{\underline p^1,\underline p^{2}}R \implies U^1U^2:=\mqty[U^1_{\mu\xi}U^2_{\xi\nu}]\in  M_{\underline p^0,\underline p^2}R.
$$

For all $\underline p\in\mathbb F$ the set of square matrices $ M_{\underline p}R:= M_{\underline p,\underline p}R$ is a $*$-ring.
The involution is given by the \textit{adjoint transpose}, i.e. $\mqty[U_{\mu\nu}]^*:=\mqty[U^*_{\nu\mu}]$.
If $R$ has a unit then  $ M_{\underline p}R$ also has one, the identity matrix $1_{\underline p}:=\mqty[\delta_{\mu\nu}]$.

For all injections $\sigma\in\mathbb F^{\text{inj}}(\underline p^0,\underline p^1)$ we define the matrix $\boldsymbol{\sigma}=\mqty[\delta_{\mu(\sigma\nu)}]\in  M_{\underline p^1,\underline p^0}\mathds Z$. We always have $\boldsymbol\sigma^*\boldsymbol\sigma=1_{\underline p^0}$. If $p^0=p^1$ then $\sigma$ is a bijection and $\boldsymbol\sigma\boldsymbol\sigma^*=1_{\underline p^1}$.

For all $p<q$ we have the natural $*$-ring monomorphism 
$$
i_{\underline p,\underline q}\in *\texttt{Ring}( M_{\underline p}\mathfrak A, M_{\underline q}\mathfrak A),
\qquad
i_{\underline p,\underline q}A:=\mqty[\begin{cases}
    A_{\mu\nu},& \max(\mu,\nu)\leq p\\
    0,&\max(\mu,\nu)>p
\end{cases}]=\mqty[A&0\\0&0].
$$
The directed colimit $ M_\infty R:=\text{colim}_\omega  M_{\underline p}R$ inherits a $*$-ring structure. The ring $ M_\infty R$ is always non-unital, even if $R$ has a unit. By identifying $R\cong  M_{\underline 1}R$ we have a natural ring inclusion $\iota:R\hookrightarrow  M_\infty R$.

The main matrix operation of relevance to $K$-theory is the direct sum.

\begin{definition}
Let $R\in *\texttt{Ring}$ and $\underline a\in\mathbb F$. The \textit{direct $\underline a$-sum} is
    \begin{align*}
        \oplus_{\underline a}: M R^{\underline a}
        \rightarrow  M R,
        \quad
        \oplus_{\underline a}U^x:=\mqty[\begin{cases}
        U^x_{\mu\nu},& x=x'\\
        0,&x\neq x'
    \end{cases}]=\mqty[U^1&0&\cdots&0\\0&U^2&\cdots&0\\
    \vdots&\vdots&\ddots&\vdots\\
    0&0&\cdots&U^x]\in M_{
        \Sigma_{\underline a}\underline p^{i0}
        ,\Sigma_{\underline a}\underline p^{i1}} R.
    \end{align*}
\end{definition}

The usual direct sum $\oplus=\oplus_{\underline 2}$ has many nice algebraic properties. It is bilinear, associative and has an identity $\oplus_{\underline 0}=0_{\underline 0}$, thus equips $MR$ with a monoid structure. This monoid structure is not commutative, but it is commutative up to similarity:
\begin{align*}
    (\boldsymbol{(12)(\underline p^{x0})})(U^1\oplus U^2)(\boldsymbol{(12)(\underline p^{x1})})^*
    &=\mqty[
        0&1_{\underline p^{20}}\\1_{\underline p^{10}}&0
    ]
    \mqty[
        U^1&0\\0&U^2
    ]
    \mqty[
        0&1_{\underline p^{11}}\\1_{\underline p^{21}}&0
    ]\\
    &=
    \mqty[
    U^2&0\\0&U^1
    ]=U^2\oplus U^1
\end{align*}

The direct sum is also compatible with compositions and adjunctions, i.e.
\begin{align*}
    (U^{11}\oplus U^{12})(U^{21}\oplus U^{22})
    &=(U^{11}U^{21})\oplus(U^{12}U^{22})\\
    (U^1\oplus U^2)^*
    &=U^{1,*}\oplus U^{2,*}
\end{align*}
whenever the dimensions of the matrices allow for matrix multiplication. We record these properties of the direct sum in the following lemma for future reference.

\begin{lemma}\label{Lema SomaComposicao}
    The direct sum is bilinear. It determines a monoid structure on $ M R$, with identity $0_{\underline 0}$, which is commutative up to similarity. The direct sum is also compatible with adjunction and composition.
\end{lemma}

Another relevant operation is the Kronecker product.

\begin{definition}
    Let $\underline m\in\mathbb F$ and $R\in *\texttt{Ring}_1$. The \textit{Kronecker $\underline m$-product} is
    $$
        \otimes_{\underline m}: MR^{\underline m}\rightarrow  M R,
        \qquad
        \otimes_{\underline m} U^i
        :=\mqty[\prod_{\underline m}U^i_{\mu^i\nu^i}]=\mqty[(\prod_{\underline{m-1}}U^i_{\mu^i\nu^i})U^m]
        \in  M_{\Pi_{\underline m}\underline p^{i0},\Pi_{\underline m}\underline p^{i1}}R.
        $$
\end{definition}

The usual Kronecker product $\otimes=\otimes_{\underline 2}$ shares some of the good algebraic properties of the direct sum. It is bilinear, associative and has $1_{\underline 1}$ as an identity element, and thus determines a monoid structure on $MR$. It also has some compatibility properties with the direct sum. The empty matrix $0_{\underline 0}$ is an absorbing element
$$
    U\otimes 0_{\underline 0}=0_{\underline 0}=0_{\underline 0}\otimes U,
$$
it strictly distributes on the right of the direct sum
$$
    (U^{11}\oplus U^{12})\otimes U^{21}
    =(U^{11}\otimes U^{21})\oplus (U^{12}\otimes U^{21}),
$$
and it distributes on the left up to similarity
\begin{align*}
    \boldsymbol{\delta_{\underline 2,\langle \underline 1,\underline 2\rangle,\langle \underline p^{ia0}\rangle}} (U^{11}\otimes (U^{21}\oplus U^{22}))\boldsymbol{\delta^*_{\underline 2,\langle \underline 1,\underline 2\rangle,\langle \underline p^{ia1}\rangle}}
    =(U^{11}\otimes U^{21})\oplus(U^{11}\otimes U^{22}).
\end{align*}

All $\delta_{\underline 2,\langle \underline 2,\underline 1\rangle,\langle \underline p^{ia}\rangle}$ are monotonous, while the $\delta_{\underline 2,\langle \underline 1,\underline 2\rangle,\langle \underline p^{ia}\rangle}$ might not be, which explains the above distinction.

Unfortunately the Kronecker product is not compatible with adjunction or composition for noncommutative $*$-rings, which greatly limits its applicability. It is also not commutative, not even up to similarity. This product does satisfy these conditions when restricted to commutative $*$-rings. The proof is as follows:
\begin{align*}
    (U^1\otimes U^2)^*
    &=\mqty[U^{2,*}_{\mu^{2}\nu^{2}}U^{1,*}_{\mu^{1}\nu^1}]\\
    &=\mqty[U^{1,*}_{\mu^{1}\nu^{1}}U^{2,*}_{\mu^{2}\nu^{2}}]\\
    &=U^{1,*}\otimes U^{2,*};\\
    (U^{11}\otimes U^{12})(U^{21}\otimes U^{22})
    &=
    \mqty[U^{11}_{\mu^1\xi^1}U^{12}_{\mu^2\xi^2}U^{21}_{\xi^1\nu^1}U^{22}_{\xi^2\nu^2}]\\
    &=\mqty[U^{11}_{\mu^1\xi^1}U^{21}_{\xi^1\nu^1}
    U^{12}_{\mu^2\xi^2}U^{22}_{\xi^2\nu^2}
    ]\\
    &=(U^{11}U^{21})\otimes(U^{12}U^{22});\\
    (\boldsymbol{(12)\langle\underline p^{i0}\rangle})(U^1\otimes U^2)(\boldsymbol{(12)\langle\underline p^{i1}\rangle})^*
    &=\mqty[U^1_{\mu^2\nu^2}U^2_{\mu^1\nu^1}]\\
    &=\mqty[U^2_{\mu^1\nu^1}U^1_{\mu^2\nu^2}]\\
    &=U^2\otimes U^1.
\end{align*}

We again record these facts in the following lemma.

\begin{lemma}\label{Prod BoasPropSempre}
    For $R\in *\texttt{Ring}_1$ the Kronecker product is bilinear, and determines a monoid structure on $ MR$ with identity $1_{\underline 1}$. The empty matrix $0_{\underline 0}$ is an absorbing element of the Kronecker product. The Kronecker product strictly distributes on the right of the direct sum and it distributes on the left up to similarity.
    
    If $R$ is commutative then the Kronecker product is commutative up to similarity, and is compatible with adjunction and composition.
\end{lemma}

\subsection{{$K$}-Theory of {$C^*$}-algebras}

The following are the main examples and constructions of $C^*$-algebras we will need to consider:

\begin{itemize}
    \item The quintessential examples of $C^*$-algebras are the algebras $\mathfrak B\mathcal H$ of bounded linear operators on Hilbert spaces. The involution is the operator adjunction 
    $A^*:=A^\dagger$ defined via the Riesz representation theorem,
    and norm the operator norm  $\lVert A\rVert_{\mathfrak B\mathcal H}
    :=\sup\limits_{\{\vec v\in \mathcal H\ \mid\ \lVert \vec v\rVert\leq 1\}}
    \lVert A\vec v\rVert_{\mathcal H}$.
    
    By the Gelfand–Naimark theorem any $C^*$-algebra $\mathfrak A$ is isometrically $*$-isomorphic to a $C^*$-subalgebra of $\mathcal B\mathcal H_{\mathfrak A}$ for some Hilbert space $\mathcal H_{\mathfrak A}$.
    
    \item Let $\mathfrak A\in \texttt C^*\texttt{Alg}_{\mathds F}$, $X\in\texttt{Top}_{\text{LCH}}$, then $C_0(X,\mathfrak A):= \texttt{Top}_*((X_+,\infty),(\mathfrak A,0))\in\texttt{C$^*$Alg}_{\mathds F}$, with $*$-algebra structure defined pointwise and norm $\lVert \varphi\rVert_{C_0(X,\mathfrak A)}:=\sup_{X}\lVert \varphi x\rVert_{\mathfrak A}$.

    For $n\in\mathds N$ the $n$-th suspension of $\mathfrak A\in \texttt{C$^*$Alg}_{\mathds F}$  is $S^n\mathfrak A:=C_0(\mathds R^n,\mathfrak A)$.
    
    \item For all $\underline p\in\mathbb F$ and $\mathfrak A\in \texttt C^*\texttt{Alg}_{\mathds C}$ the matrix $*$-algebra $ M_{\underline p}\mathfrak A$ is a  $C^*$-algebras, with norm the operator norm
$$
    \left\lVert A\right\rVert_{ M_{\underline p}\mathfrak A}
        :=\left\lVert  A\right\rVert_{\mathcal B(\mathcal H^{\underline p}_{\mathfrak A})}.
$$ 

\item Let $\mathfrak A\in \texttt C^*\texttt{Alg}_{\mathds F}$. The ring of finite matrices $ M_\infty \mathfrak A$ is normed. The \textit{stabilization} $\mathfrak K \mathfrak A$ of $\mathfrak A$ is the norm-completion of $ M_\infty \mathfrak A$. We will denote elements of $\mathfrak{KA}$ with a hat, i.e. $\hat A\in\mathfrak{KA}$.

\item For any $C^*$-algebra $\mathfrak A$ the unital $*$-algebra $\widetilde{\mathfrak A}$ equipped with the norm 
$$\lVert \alpha\tilde 1+A\rVert_{\widetilde{\mathfrak A}}:=\lVert \alpha \text{id}+L_{A}\rVert_{\mathcal B\mathfrak A},
$$
where $L_{A}B=AB$, is a $C^*$-algebra.

For all $\mathfrak A\in \texttt C^*\texttt{Alg}_{\mathds F}$ we have a split short exact sequence of $C^*$-algebras
\begin{gather*}
    \xymatrix{\mathfrak A\ar@{^{(}->}[r]^\iota&\widetilde{\mathfrak A}\ar@<0.1cm>@{->>}[r]^{\pi}&\mathds F\ar@<0.1cm>[l]^\lambda};\\
    \iota A:=0\tilde 1+A,
    \qquad
    \pi_{\mathds C}(\alpha \tilde 1+A):=\alpha,
    \qquad
    \lambda \alpha:=\alpha\tilde 1+0.
\end{gather*}

The \textit{scalar map} $s\in \texttt C^*\texttt{Alg}_{\mathds F}(\widetilde{\mathfrak A},\widetilde{\mathfrak A})$ is defined as $s:=\lambda\pi$, so $s(\alpha \tilde 1+A)=\alpha\tilde 1+0$.

For $\mathfrak A\in \texttt C^*\texttt{Alg}_1$ and $\underline p\in\mathbb F$ we have distinct matrices $1_{\underline p},\tilde 1_{\underline p}\in  M_{\underline p}\widetilde{\mathfrak A}$ with 
$$
1_{\underline p}=\mqty[0\tilde 1+\delta_{\mu'\mu}]
,\qquad
\tilde 1_{\underline p}=\mqty[\delta_{\mu'\mu}\tilde 1+0],
$$the latter being the ring unit.

\item We will be particularly interested in matrices in $ M\widetilde{\mathfrak{KA}}$, which we will denote using bold type: 
$$
    \boldsymbol A\equiv\left[ \alpha_{\mu\nu} \tilde 1 +\hat A_{\mu\nu}\right]\in  M_{\underline p^1,\underline p^0} \widetilde{\mathfrak{KA}}
$$

\end{itemize}

In operator $K$-theory we study unital $C^*$-algebras $\mathfrak{A}$ by analyzing projection matrices, i.e. self-adjoint idempotent matrices. For all $\underline p\in\mathbb F$ and $\mathfrak A\in \texttt C^*\texttt{Alg}_{\mathds F}$ the space of projection matrices of order $\underline p$ over $\mathfrak A$ is 
$$
    \mathcal P_{\underline p}\mathfrak A:=\{P\in  M_{\underline p}\mathfrak A\mid P=P^*=P^2\}.
$$
We then define 
$$\textstyle
    \mathcal P_\infty\mathfrak A:=\coprod_{\mathbb F}\mathcal P_{\underline p}\mathfrak A.
$$

The range of each $P\in \mathcal P_\infty\mathfrak A$ is a finitely generated projective right $\mathfrak A$-module, and up to isomorphism any such module is isomorphic to the image of a projection in $\mathcal P_\infty\mathfrak A$. If $\mathfrak A$ is commutative then for each $P\in \mathcal P_\infty\mathfrak A$ Swan's theorem tells us that 
$\Sigma_{\text{Sp}\mathfrak A}\chi P$ is a finite dimensional vector bundles over the Gelfand spectrum of characters $\text{Sp}\mathfrak A$.

By lemma \ref{Lema SomaComposicao} the direct sum preserves projections and induce a monoid structure on $(\mathcal P_\infty\mathfrak A;\oplus,0_{\underline 0})$. This monoid structure is not commutative, but it does induce a commutative monoid structure on the  set of Murray-von Neumann equivalence classes. The Murray-von Neumann relation $\sim_0$ is defined on $\mathcal P_\infty\mathfrak A$ as follows. If $P,Q\in \mathcal P_\infty \mathfrak A$ then $P\sim_0 Q$ if there is a partial isometriy\footnote{Partial isometries are elements $U\in\mathfrak A$ such that $U^*U\in\mathcal P_\infty\mathfrak A$, or equivalently such that $UU^*\in\mathcal P_\infty\mathfrak A$.} $U\in M_{\underline q,\underline p}\mathfrak A$ satisfying 
$$
    U^*U=P, \qquad UU^*=Q,
$$
so in particular $\text{Im}P=\text{Ker}U^\perp$ and $\text{Im}U=\text{Im}Q$. 

We denote the Murray-von Neumann equivalence class of a projection matrix $P$ as $[P]_0$, and the set of Murray-von Neumann equivalence classes as $\widetilde{\texttt{pr}}_{\mathfrak A}:=\mathcal P_\infty\mathfrak A_{/\sim_0}$. For $P^1,P^2\in \mathcal P_\infty \mathfrak A$ we have that
$$
\boldsymbol{(12)(\underline p^x)}(P^1\oplus P^2)=\mqty[0&P^2\\P^1&0]
$$
is a partial isometry witnessing $P^1\oplus P^2\sim_0 P^2\oplus P^1$. Therefore the operation 
$$[P^1]_0+[P^2]_0:=[P^1\oplus P^2]_0$$
determines a commutative monoid structure on $\widetilde{\texttt{pr}}_{\mathfrak A}$.

It will be useful to organize the structure of the projection matrices, partial isometries witnessing the Murray-von Neumann relations between them and the direct sum into a natural permutative category construction.

\begin{definition}
    The \textit{Murray-von Neumann permutative category functor} is
    \begin{gather*}
        \texttt{pr}_{-}:\texttt{C$^*$Alg}_{\mathds F}\rightarrow \texttt{PermCat};\\
        \text{Ob}\ \texttt{pr}_{\mathfrak A}:=\mathcal P_\infty\mathfrak A,
        \qquad \texttt{pr}_{\mathfrak A}(P,Q)
        :=\left\{U\in M_{\underline q,\underline p}\mathfrak A
        \mid
        P=U^*U,\ UU^*=Q\right\},\\
        \text{id}_P:=P, 
        \qquad U^1U^2:=\mqty[U^1_{\mu\xi}U^2_{\xi\nu}];\\
        U^1\oplus U^2:=\mqty[U^1&0\\0&U^2],
        \quad
        \boldsymbol 0:=0_{\underline 0},
        \quad 
        \tau_{P^1,P^2}:=\boldsymbol{(12)(\underline p^x)}P^1\oplus P^2=\mqty[0&P^2\\P^1&0].
    \end{gather*}
        
\end{definition}

By construction we have a natural isomorphism $\widetilde{\texttt{pr}}_{\mathfrak A}\cong \pi_0\lvert \texttt{pr}_{\mathfrak A}\rvert$.

\begin{definition}
    The \textit{$K_{00}$-group functor} is
    \begin{gather*}
        K_{00}:\texttt C^*\texttt{Alg}_{\mathds F,1}\rightarrow \texttt{AbGrp},
        \qquad
        K_{00} \mathfrak A:= \text{Gr}\ \widetilde{\texttt{pr}}_{\mathfrak A}.
    \end{gather*}
\end{definition}

This functor is half-exact.
The definition of $K_{00}$ could be made on the whole of $\texttt C^*\texttt{Alg}_{\mathds F}$, but then it wouldn't have this important homological property. 

\begin{definition}
    The \textit{graded $K$-group functor} is
$$
    K_-:\texttt C^*\texttt{Alg}_{\mathds F}\rightarrow \texttt{AbGrp}^{\mathds N},
    \qquad K_n\mathfrak A:=\begin{cases}
        \text{Ker}(K_{00}\pi_{\mathds C}:K_{00} \widetilde{\mathfrak A} 
        \rightarrow 
        K_{00}\mathds C),& n=0\\
        K_0 S^n\mathfrak A, &n>0
    \end{cases}.
$$
\end{definition}

The functor $K_-$ is homotopy invariant, i.e. if $H\in\texttt{C$^*$Alg}_{\mathds F}(\mathfrak A,C(I,\mathfrak B))$ then 
        $$(\text{ev}_0H)_*=(\text{ev}_1 H)_*\in\texttt{AbGrp}^{\mathds N}(K_-\mathfrak A,K_-\mathfrak B).
        $$
For any short exact sequences
    $\xymatrix{
    \mathfrak I
    \ar@{^{(}->}[r]^{\iota}
    &\mathfrak A
    \ar@{->>}[r]^{p}
    &\mathfrak A_{/\mathfrak I}
}$ there is a natural boundary transformation 
$\partial:K_{n+1}\mathfrak A_{/\mathfrak I}\Rightarrow K_n\mathfrak I$
such that we get a natural long exact sequence
$$
    \xymatrix{\cdots\ar[r]
    &K_{n+1}\mathfrak A_{/\mathfrak I}
    \ar[r]^{\partial}
    &K_{n}\mathfrak I
    \ar[r]^{\iota_*}
    &K_{n}\mathfrak A
    \ar[r]^{p_*}
    &K_{n}\mathfrak A_{/\mathfrak I}
    \ar[r]&\cdots}.
$$
Thus $K_-$ forms an extraordinary homology theory. The functor $K_-$ preserves directed colimits. In particular this implies that the natural inclusion $\iota:\mathfrak A\hookrightarrow \mathfrak K \mathfrak A$ induces an isomorphism $\iota_*:K_- \mathfrak A\cong K_-  \mathfrak K \mathfrak A$.

An important property of the $K$-groups is Bott periodicity. For complex $C^*$-algebras we have natural isomorphisms $K_n\mathfrak A\cong K_{n+2}\mathfrak A$, and for real $C^*$-algebras we have natural isomorphisms $K_n\mathfrak A\cong K_{n+8}\mathfrak A$. This is a useful property since this means we need to compute only a finite number of groups to understand the $K$-theory of $C^*$-algebras, and the long exact sequences induced by short exact sequences collapse to finite exact cycles. Bott periodicity also allows us to extend $K_-\mathfrak A$ to a $\mathds Z$-graded abelian group, with negative $K$-groups inductively defined. We will return to this issue in the last subsection of this article.

If $\mathfrak A\in \texttt C^*\texttt{Alg}$ is commutative then the Kronecker product induces a \textit{cup product}
\begin{gather*}
    \cup:K_{n^1}\mathfrak A\otimes K_{n^2}\mathfrak A\rightarrow K_{n^1+n^2}\mathfrak A
\end{gather*}
that equips $K_-\mathfrak A$ with a graded commutative ring structure, which is associated by Gelfand duality and Swan's theorem to the graded ring structure of topological $K$-theory \cite{karoubi2008k}. It is this structure we wish to generalize to the non-commutative context.

\section{The $\mathscr L$-permutative categories {$\texttt{pr}_{\widetilde{\mathfrak{KA}}}$}}

In this section we define our generalization of bipermutative categories, which will take the form of actions of operad pairs on categories. We first review basic definitions of operads, operad pairs and their algebras \cite{MaEinftyRingSpcsSpectra,may2006geometry,Ma09a}. We then define $\mathscr G$-categories and $(\mathscr P,\mathscr G)$-categories in a similar way, with equivariance and distributivity satisfied only up to natural isomorphism. We then define $\mathscr L$-permutative categories as $(\mathscr C,\mathscr L)$-categories and prove that $\texttt{pr}_{\widetilde{\mathfrak{KA}}}$ have natural $\mathscr L$-permutative category structures.

\subsection{Operads and their algebras}

\begin{definition}
    A \textit{topological operad} is a contravariant functor 
    $$\mathscr G:\mathbb F^{\text{bij,op}}\rightarrow \texttt{Top}
    $$
    equipped with an abstract identity element and a composition map
$$
\textbf{id}\in\mathscr G\underline 1, 
\qquad \circ=\langle\circ_{\underline m,\langle \underline n^i\rangle}\rangle\in  \prod_{\Sigma_{\mathbb F}\mathbb F^m}\texttt{Top}\left(\mathscr G\underline m\times\prod_{\underline m}\mathscr G\underline n^i,\mathscr G\Sigma_{\underline m}\underline n^i\right)
$$
such that $\lvert \mathscr G\underline 0\rvert=1$ and, using the notation $f\langle g^i\rangle:=\circ(f,\langle g^i\rangle)$, the following associativity, unit and equivariance conditions are satisfied:
\begin{align*}
    f\langle g^i \langle h^{ij}\rangle\rangle&=(f \langle g^i\rangle)\langle h^{ij}\rangle;\\
    \textbf{id} f=&\ f=f\langle \textbf{id}\rangle;\\
    (f\cdot\sigma)\langle g^i\rangle
    &=(f\langle g^{\sigma^{-1}i}\rangle)\cdot\sigma(\underline n^i);\\
    f\langle g^i\cdot\sigma^i\rangle
    &=(f\langle g^i\rangle)\cdot\Sigma_{\underline m}\sigma^i.
\end{align*}
\end{definition} 

For each $X\in \texttt{Top}$ we have the functor
$$
X^-:\mathbb F^{\text{bij}}\rightarrow \texttt{Top}, \qquad \sigma\cdot\langle x^i\rangle:=\langle x^{\sigma^{-1}i}\rangle.
$$
For an  operad $\mathscr G$ we define the monad $(G;\eta,\eta)$ on $\texttt{Top}$ as
\begin{gather*}
    GX:=\int^{\mathbb F^{\text{bij}}}\mathscr G\underline m\times X^{\underline m};
    \qquad \eta x:=[\textbf{id},x],
    \quad \mu[f,[g^i,\langle x^{ij}\rangle]]:=[f\langle g^i\rangle, \langle x^{ij}\rangle].
\end{gather*}

\begin{definition}
    A \textit{$\mathscr G$-space} is a $G$-algebra. The category of $\mathscr G$-spaces is denoted $\mathscr G[\texttt{Top}]$. Equivalently a $\mathscr G$-space is a topological space $X$ equipped with a $\mathscr G$-action, i.e. a map
$$
\xi=\langle \xi_{\underline m}\rangle\in \prod_{\mathbb F}\texttt{Top}\left(\mathscr G\underline m\times X^{\underline m}, X\right)
$$
satisfying, using the notation $f\langle x^i\rangle:=\xi(f,\langle x^i\rangle)$, the following equivariant action conditions:
\begin{align*}
f\langle g^i\langle x^{ij}\rangle\rangle
&=(f\langle g^i\rangle)\langle x^{ij}\rangle;\\
\textbf{id} x&=x;\\
(f\cdot\sigma)\langle x^i\rangle&=f( \sigma\cdot\langle x^i\rangle).
\end{align*}
\end{definition}

The underlying functor of an operad $\mathscr G$ can be extended to a functor on $\mathbb F^{\text{inj},\text{op}}$. For $\sigma\in\mathbb F^{\text{inj}}(\underline m^0,\underline m^1)$ the right action $\cdot\sigma\in\texttt{Top}(\mathscr G\underline m^1,\mathscr G\underline m^0)$ is
$$
f\cdot\sigma:=f\left\langle\begin{cases}
e,&i^1\not\in\text{Im }\sigma\\
\textbf{id},&i^1\in\text{Im }\sigma
\end{cases}\right\rangle.
$$
These morphisms are the \textit{degenerations} of the operad. We denote by $e_X$ the element in $X\in \mathscr G[\texttt{Top}]$ picked out by $e\in\mathscr G\underline 0$. For each $(X,e_X)\in \texttt{Top}_*$ we have the functor
$$
X^-:\mathbb F^{\text{inj}}\rightarrow \texttt{Top}_*, \qquad \sigma\cdot\langle x^i\rangle=\left\langle\begin{cases}
    e_X,&i'\not\in \text{Im }\sigma\\
    x^{\sigma^{-1}i'},&i'\in \text{Im }\sigma
\end{cases}\right\rangle.
$$

We denote by $\mathscr C$ the operad with a unique element in each $\mathscr C \underline m$, which by abuse of notation we will denote simply as $\underline m\in \mathscr C \underline m$. It is the terminal object in the category of operads, and the $\mathscr C$-spaces are precisely the topological commutative monoids.

\begin{definition}
    An \textit{$E_\infty$-operad} $\mathscr E$ is an operad such that each space $\mathscr E\underline m$ is a $\mathbb F^{\text{bij}}(\underline m,\underline m)$-free contractible space.
\end{definition}

An important example of $E_\infty$-operad is the linear isometries operad $\mathscr L$, defined as
\begin{gather*}
    \mathscr L\underline m
    :=\{f=\langle f_i\rangle\in\mathscr I(\oplus_{\underline m}\mathds R^\infty,\mathds R^\infty)\};\\
    f\cdot\sigma:=\langle f_{\sigma i}\rangle,
    \qquad\textbf{id}:=\text{id}_{\mathds R^\infty},
    \qquad f\langle  g^i\rangle:=\langle f_ig^i_j\rangle.
\end{gather*}

\begin{definition}
    An \textit{operad pair} $(\mathscr P,\mathscr G,\ltimes)$ is a pair of operads equipped with an \textit{action} of $\mathscr G$ on $\mathscr P$, i.e. a map
$$
\ltimes=\langle \ltimes_{\underline m,\langle \underline a\rangle}\rangle \in \prod_{\Sigma_{\mathbb F}\mathbb F^m}\texttt{Top}\left(\mathscr G \underline m \times\prod_{\underline m}\mathscr P \underline a^i,\mathscr P \Pi_{\underline m}\underline a^i\right)
$$
such that, using the notations $f\ltimes \langle \alpha^i \rangle:=\ltimes(f,\langle \alpha^i\rangle)$, $0\in \mathscr P\underline 0$ and $1\in \mathscr G\underline 0$, the following distributivity, unit and equivariance conditions:
\begin{align*}
    f\ltimes\langle g^i\ltimes\langle \alpha^{ij}\rangle\rangle
    &=(f\langle g^i\rangle)\ltimes\langle \alpha^{ij}\rangle; \\
    f\ltimes \langle \alpha^i\langle\beta^{ix}\rangle\rangle
    &=((f\ltimes\langle\alpha^i\rangle)\langle f\ltimes\langle\beta^{ix^i}\rangle\rangle)\cdot \delta_{\underline m,\langle \underline a^i\rangle,\langle\underline b^{ix}\rangle}; \\
    \textbf{id}\ltimes \alpha
    &=\alpha;\\
    f\ltimes\langle \textbf{id}\rangle
    &=\textbf{id};\\
    (f\cdot\sigma)\ltimes\langle\alpha^{i}\rangle
    &=(f\ltimes\langle\alpha^{\sigma^{-1}i}\rangle)\cdot\sigma\langle \underline a^{i}\rangle;\\
    f\ltimes\langle\alpha^i\cdot\sigma^i\rangle
    &=(f\ltimes\langle \alpha^i\rangle)\cdot\Pi_{\underline m}\sigma^i;\\
    1\ltimes \ast
    &=\textbf{id};\\
    \exists i\in\underline m:\alpha^i=0\implies f\ltimes\langle\alpha^i\rangle
    &=0.
\end{align*}
The operad $\mathscr P$ is the \textit{additive operad}, and $\mathscr G$ is the \textit{multiplicative operad}.
\end{definition}

All operads $\mathscr G$ act on $\mathscr C$, with 
$$f\ltimes\langle \underline a^i\rangle:=\Pi_{\underline m}\underline a^i.
$$

Given an operad pair $(\mathscr P,\mathscr G,\ltimes)$ the monad $P$ on $\texttt{Top}$ induces a monad, also denoted $P$, on $\mathscr G[\texttt{Top}]$. 

\begin{definition}
    A \textit{$(\mathscr P,\mathscr G)$-space} is a $P$-algebra in $G[\texttt{Top}]$.
    
    Equivalently a $(\mathscr P,\mathscr G)$-space is a topological space $Z$ equipped with a $\mathscr P$-space and a $\mathscr G$-space structures such that the $\mathscr P$-neutral element $0_Z$ is an absorbing element for the $\mathscr G$-space structure, i.e.
$$
\exists i\in\underline m:z^i=0\implies f\langle z^i\rangle=0_Z,
$$
and such that the following distributivity equation holds:
$$
f\langle\alpha^i\langle z^{ix}\rangle\rangle
=(f\ltimes\langle\alpha^i\rangle)\langle f\langle z^{ix^i}\rangle\rangle.
$$
\end{definition}

\begin{definition}
    An \textit{$E_\infty$-ring space} is a $(\mathscr P,\mathscr G)$-space for an operad pair $(\mathscr P,\mathscr G,\ltimes)$ where both $\mathscr P$ and $\mathscr G$ are $E_\infty$-spaces.
\end{definition}

We will be particularly interested in $(\mathscr H_\infty,\mathscr L)$-spaces, where $\mathscr H_\infty$ is Steiner's operad of paths of embeddings \cite{St79}. For $U\subset \mathscr A$ the operad of $U$-embeddings $\text{Emb}_U$ is
\begin{gather*}
    \text{Emb}_U\underline a:=\{\alpha=\langle\alpha_x\rangle\in U^{\sqcup_{\underline a} U}\mid \alpha\text{ is an embedding}\};\\
    \alpha \cdot\sigma:=\langle \alpha_{\sigma x}\rangle,
    \qquad  \textbf{id}:=\text{id}_U,
    \qquad \alpha\langle \beta^x\rangle:=\langle \alpha_x\beta^x_y\rangle.
\end{gather*}
For all $U\subset \mathscr A$ the $U-$Steiner operad $\mathscr H_U$ is
\begin{gather*}
    \mathscr H_U\underline a:= \left\{\alpha=\langle\alpha_x\rangle\in U^{I\times \sqcup_{\underline a} U}\left\lvert\begin{matrix*}[l]
    \forall x\in \underline a,\forall t\in I,\forall \vec u,\vec v\in U:\\
    (\vec w\mapsto\alpha_x(t,\vec w))\in \text{Emb}_U\{x\},\\
    \lVert \alpha_x(t,\vec u)- \alpha_a(t,\vec v)\rVert\leq \lVert\vec u- \vec v\rVert,\\
    \alpha_x(1,\vec u)=\vec u,\\
    \langle \vec w\mapsto\alpha_{a'}(0,\vec w)\rangle\in \text{Emb}_U \underline a.
    \end{matrix*}\right.\right\}\\
    \alpha \cdot\sigma:=\langle\alpha_{\sigma x}\rangle, 
    \qquad\textbf{id}:=\text{pr}_U,
    \qquad \alpha \langle \beta^a\rangle:=\langle  (t,\vec u)
    \mapsto\alpha_x(t,\beta^x_y(t,\vec u))\rangle
\end{gather*}
For all $(U,V)\in\Sigma_{\mathscr A}\mathscr A_U$ we have a natural inclusions
$$
    \iota_{U,V}:\mathscr H_U\Rightarrow\mathscr H_V,
    \quad\iota_{U,V} \alpha :=\langle t\mapsto \text{pr}_{V-U}+(\alpha_xt)\text{pr}_U\rangle,
$$
and Steiner's operad is defined as $\mathscr H_\infty:=\text{colim}_{U\subset \mathds R^\infty}\mathscr H_U$. The action of $\mathscr L$ on $\mathscr H_\infty$ is
$$\textstyle
    f\ltimes\langle \alpha^i\rangle:=
    \left\langle (t,\vec u)\mapsto \text{pr}_{f^\perp}+\sum_{\underline m} f_i \alpha^{i}_{x^i}(t,\vec u) f_i^* \right\rangle.
$$

\subsection{$(\mathscr P,\mathscr G)$-categories}

We introduce the following characterization of categorial algebras over a topological operad.

\begin{definition}
    Let $\mathscr G$ be a topological operad. A \textit{$\mathscr G$-category} $\mathcal C$ consists of
    \begin{itemize}
        \item A topological category $\mathcal C$;
        
        \item A sequence of functors 
        $$
            \circledast=\langle\circledast_{\underline m}\rangle\in \prod_{\mathbb F} \texttt{Cat}\left(\mathscr G\underline m\times\mathcal C^{\underline m}, \mathcal C\right)
        $$
        where we consider each $\mathscr G\underline m$ as a trivial topological category, i.e. containing only identity morphisms. We will use the notation $\circledast_f:=\circledast_{\underline m}(f,-)$.
        
        \item For each $\sigma\in\Sigma_{\mathbb F^2}\mathbb F^{\text{bij}}(\underline m,\underline m)$ a natural isomorphism
        $$\xymatrix@=2cm{
            \mathscr G\underline m\times\mathcal C^{\underline m}
            \rtwocell^{\circledast_{f\cdot\sigma}P^{i}}_{\circledast_{f}P^{\sigma^{-1}i}}{\quad\ \tau^{\circledast,\sigma}_{f,\langle P^i\rangle}}
            &\mathcal C.
        }$$
    \end{itemize}
    
    Such that the following cohesion equations hold:
    \begin{align*}
        \circledast_{f\langle g^i\rangle}
        &\textstyle=\circledast_f\langle\circledast_{g^i}\rangle;\\
        \circledast_{\textbf{id}}&=\text{id}_{\mathcal C};\\ \tau^{\circledast,\sigma}_{f,\langle P^{\sigma^{\prime-1}i}\rangle}\tau^{\circledast,\sigma'}_{f\cdot\sigma,\langle P^i\rangle}&=\tau^{\circledast,\sigma\sigma'}_{f,\langle P^i\rangle};\\
        \tau^{\circledast,\text{id}}_{\alpha,\langle P^i\rangle}&=\text{id}_{\circledast_\alpha P^i}.
    \end{align*}
    The third equation equivalent to the commutation of the natural cohesion diagram
    $$
    \xymatrix{\circledast_{f\cdot\sigma\sigma'}P^i
    \ar@{=>}[d]_{\tau^{\circledast,\sigma'}_{f\cdot\sigma,\langle P^i\rangle}}
    \ar@{=>}[dr]^{\tau^{\circledast,\sigma\sigma'}_{f,\langle P^i\rangle}}&\\
    \circledast_{f\cdot\sigma}P^{\sigma^{\prime-1}i}
    \ar@{=>}[r]_{\tau^{\circledast,\sigma}_{f,\langle P^{\sigma^{\prime-1}i}\rangle}}&\circledast_f P^{(\sigma\sigma')^{-1}i}}.
    $$
    
    We denote the category of $\mathscr G$-categories as $\mathscr G[\texttt{Cat}]$.
\end{definition}

\begin{proposition}\label{Prop isom Perm CCat}
    The categories $\texttt{PermCat}$ and $\mathscr C [\texttt{Cat}]$ are isomorphic.
\end{proposition}

\textbf{Proof:} For $(\mathcal C;\circledast,\boldsymbol e,\tau^{\circledast}_{P^1,P^2})\in \texttt{PermCat}$ we can define $\circledast_{\underline 0}:=\boldsymbol e$, $\circledast_{\underline 1}:=\text{id}_{\mathcal C}$, $\circledast_{\underline 2}:=\circledast$ and for each $m>2$ we inductively define $\circledast_{\underline m}:=\circledast (\circledast_{\underline{m-1}}
\times \text{id}_{\mathcal C})$. The coherence conditions on the braiding transformation $\tau$ means there is a unique choice for the $\tau^{\circledast,\sigma}_{\underline a,\langle P^x\rangle}$. The strict associativity and identity conditions and the coherence of the braiding implies these define a $\mathscr C $-category.

For $(\mathcal C;\circledast_{\underline a},\tau^{\circledast,\sigma}_{\underline a,\langle P^x\rangle})\in \mathscr C [\texttt{Cat}]$ we can define $\circledast:=\circledast_{\underline 2}$, $\boldsymbol e:=\circledast_{\underline 0}$ and $\tau^{\circledast}_{P^1,P^2}:=\tau^{\circledast,(12)}_{\underline 2,\langle P^x\rangle}$. The operad action conditions implies this defines a permutative category.$\blacksquare$\\

\begin{definition}
    If $(\mathscr P,\mathscr G,\ltimes)$ is an operad pair then a \textit{$(\mathscr P,\mathscr G)$-category} consists of a topological category $\mathcal C$ equipped with
    \begin{itemize}
        \item A $\mathscr P$-category structure $(\mathcal C;\oplus_{\alpha},\tau^{\oplus,\sigma}_{\alpha,\langle P^x\rangle})$
        
        \item A $\mathscr G$-category structure $(\mathcal C;\otimes_f,\tau^{\otimes,\sigma}_{f,\langle P^i\rangle})$
        
        \item For each $(\underline m,\langle\underline a^i\rangle)\in\Sigma_{\mathbb F}\mathbb F^m$ a distributivity natural transformation
    $$\xymatrix@=2cm{
        \mathscr G\underline m
        \times \prod_{\underline m}\mathscr P\underline a^i
        \times\mathcal C^{\Sigma_{\underline m}\underline a^i}
        \rtwocell^{\qquad\otimes_f\oplus_{\alpha^i}P^{ix}}_{\qquad\oplus_{f\ltimes \langle \alpha^i\rangle}\otimes_f P^{ix^i}}{\qquad\ \delta_{f,\langle\alpha^i \rangle,\langle P^{ix}\rangle}}
        &\mathcal C
    }$$
    \end{itemize}
    such that $\oplus_0$ is an absorbing element for $\otimes_-$, i.e.
    \begin{gather*}
        \exists i\in\underline m:
        \alpha^i= 0
        \implies
        \otimes_f \oplus_{\alpha^i}=\oplus_0,
    \end{gather*}
    and the following coherence diagrams commute:
    \begin{gather*}
    \xymatrix@C=1.75cm{
        \otimes_f\oplus_{\alpha^i\cdot\sigma^i}  P^{ix}
        \ar@{=>}[d]_{\otimes_f \tau^{\oplus,\sigma^i}_{\alpha^i,\langle  P^{ix}\rangle}}
        \ar@{=>}[r]^{\delta_{f,\langle \alpha^i\cdot\sigma^i\rangle,\langle  P^{ix}\rangle}}
        &\oplus_{f\ltimes\langle\alpha^i\cdot\sigma^i\rangle}\otimes_f  P^{ix^i}
        \ar@{=>}[d]^{\tau^{\oplus,\Pi_{\underline m}\sigma^i}_{f\ltimes\langle\alpha^i\rangle,\langle \otimes_f P^{ix^i}\rangle}}\\
        \otimes_f\oplus_{\alpha^i} P^{i(\sigma^{i,-1}x)}
        \ar@{=>}[r]_{\delta_{f,\langle \alpha^i\rangle,\langle  P^{i\ \sigma^{i,-1}x}\rangle}}
        &\oplus_{f\ltimes\langle\alpha^i\rangle}\otimes_f  P^{i(\sigma^{i,-1}x^i)}},
        \quad\xymatrix@C=1.75cm{\otimes_{f\cdot\sigma}\oplus_{\alpha^{i}}  P^{ix}
        \ar@{=>}[r]^{\delta_{f\cdot\sigma,\langle \alpha^{i}\rangle,\langle  P^{ix}\rangle}}
        \ar@{=>}[dd]_{\tau^{\otimes,\sigma}_{f,\langle \oplus_{\alpha^{i}} P^{ix}\rangle}}
        &\oplus_{f\cdot\sigma\ltimes\langle\alpha^{i}\rangle}\otimes_{f\cdot \sigma}  P^{ix^{i}}
        \ar@{=>}[d]^{\oplus_{f\cdot\sigma\ltimes\langle\alpha^{i}\rangle}\tau^{\otimes,\sigma}_{f,\langle  P^{ix^{i}}\rangle}}\\
        &\oplus_{f\cdot\sigma\ltimes\langle\alpha^{i}\rangle}\otimes_{f}  P^{(\sigma^{-1}i)x^{\sigma^{-1}i}}
        \ar@{=>}[d]^{\tau^{\oplus,\sigma\langle \underline a^{i}\rangle}_{f\ltimes\langle\alpha^{\sigma^{-1}i}\rangle,\langle \otimes_{f} P^{(\sigma^{-1}i)x^{\sigma^{-1}i}}\rangle}}\\
        \otimes_{f}\oplus_{\alpha^{\sigma^{-1}i}}  P^{(\sigma^{-1}i)x}
        \ar@{=>}[r]_(0.45){{\delta_{f,\langle \alpha^{\sigma^{-1}i}\rangle,\langle  P^{(\sigma^{-1}i)x}\rangle}}}
        &\oplus_{f\ltimes\langle\alpha^{\sigma^{-1}i}\rangle}\otimes_{f}  P^{(\sigma^{-1}i)x^{i}}},\\
        \xymatrix@C=1.75cm{
        \otimes_f\oplus_{\alpha^i}\oplus_{\beta^{ix}}  P^{ixy}
        \ar@{=}[d]
        \ar@{=>}[r]^(0.45){\delta_{f,\langle  \alpha^i\rangle,\langle \oplus_{\beta^{ix}} P^{ixy}\rangle}}&
        \oplus_{f\ltimes\langle\alpha^i\rangle}\otimes _f\oplus_{\beta^{ix^i}}  P^{ix^iy}
        \ar@{=>}[d]^{\oplus_{f\ltimes\langle\alpha^i\rangle}\delta_{f,\langle  \beta^{ix^i}\rangle,\langle  P^{ix^iy}\rangle}}\\
        \otimes_f\oplus_{\alpha^{i}\langle \beta^{ix}\rangle}  P^{ixy}
        \ar@{=>}[d]_{\delta_{f,\langle  \alpha^i\langle \beta^{ix}\rangle\rangle,\langle  P^{ixy}\rangle}}&
        \oplus_{f\ltimes\langle\alpha^i\rangle}\oplus_{f\ltimes\langle \beta^{ix^i}\rangle}\otimes _f  P^{ix^iy^i}
        \ar@{=}[d]\\
        \oplus_{f\ltimes\langle\alpha\langle\beta^{ix}\rangle\rangle}\otimes_f  P^{ix^iy^i}
        \ar@{=>}[r]_(0.55
        ){\tau^{\oplus,\delta_{\underline m,\langle  \underline a^i\rangle,\langle \underline b^{ix}\rangle}}_{f\ltimes\langle\alpha^i\rangle \langle f\ltimes\langle \beta^{ix^i}\rangle\rangle,\langle \otimes_f  P^{ix^iy^i}\rangle}}
        &\oplus_{f\ltimes\langle\alpha^i\rangle \langle f\ltimes\langle \beta^{ix^i}\rangle\rangle}\otimes _f  P^{ix^iy^i}
    }.
    \end{gather*}
\end{definition}

These coherence conditions were designed after the ones for bipermutative categories, and we can construct an isomorphism $(\mathscr C,\mathscr C)[\texttt{Cat}]\cong \texttt{BipCat}$ similar to the one in the proof of proposition \ref{Prop isom Perm CCat}.

\begin{definition}
    An $\mathscr L$-\textit{permutative category} is a $(\mathscr C,\mathscr L)$-category.
\end{definition}

By pullback on the terminal operad morphism $\mathscr L\rightarrow \mathscr C$ all bipermutative categories are $\mathscr L$-permutative categories. The following $\mathscr L$-permutative categories will be relevant to our main results.

\begin{definition}
    For $\mathds F\in \{\mathds R,\mathds C\}$ the $\mathds F$-\textit{unitary bipermutative category} $\mathcal U_{\mathds F}$ is defined as
\begin{align*}
    \noalign{\centering$\text{Ob}\ \mathcal U_{\mathds F}:=\mathbb F, \qquad \mathcal U_{\mathds F}(\underline p,\underline q):=\left\{\begin{cases}
        \{U\in M_{\underline p}\mathds F\mid 1_{\underline p}=U^*U=UU^*\},&p=q,\\
        \emptyset , & p\neq q.
    \end{cases}\right\};$}
    \oplus_{\underline a} U^x
    &:=\left[ \begin{cases}
        U^x_{\mu\nu},&x=x'\\
        0,& x\neq x'
    \end{cases}
    \right]
    ,&
    \tau^{\oplus,\sigma}_{\underline a,\langle \underline p^x\rangle}
    &:=\boldsymbol{\sigma (\underline p^x)};\\
    \otimes_{\underline m} U^i
    &:=\mqty[\prod_{\underline m}U^{i}_{\mu^i\nu^i}]
    ,&
    \tau^{\otimes,\sigma}_{\underline m,\langle \underline p^i\rangle}
    &:=\boldsymbol{\sigma \langle\underline p^i\rangle};\\
    \noalign{\centering$\delta_{\underline m,\langle \underline a^i\rangle,\langle \underline p^{ix}\rangle}:=\boldsymbol{\delta_{\underline m,\langle \underline a^i\rangle,\langle \underline p^{ix}\rangle}}.$}
    \end{align*}
\end{definition}

If $\mathfrak A\in \texttt{C$^*$Alg}_{\mathds F}$ is commutative and unital then there is a bipermutative structure on $\texttt{pr}_{\mathfrak A}$, with multiplication defined by the Kronecker product. In general, even though $ \texttt{pr}_{\mathfrak A}$ is not closed under the Kronecker product, we can equip the unitization of its stabilization with an $\mathscr L$-permutative category structure. 

\begin{lemma}
    The monoid $\mathscr L\underline 1$ acts on $ \widetilde{\mathfrak{KA}}$ by 
    $f\cdot(\upsilon\tilde 1+\hat U):= \upsilon\tilde 1 + f\hat Uf^*$.
    
    For all $ (f,\langle \upsilon^i\tilde 1+\hat U^i\rangle)\in \mathscr L\underline 2\times \widetilde{\mathfrak{KA}}^{\underline 2}$ we have
\begin{align*}
    (\upsilon^1\tilde 1+f_1\hat U^1f_1^*)(\upsilon^2\tilde 1+f_2\hat U^2f_2^*)&=
    \upsilon^1\upsilon^2\tilde 1 + \upsilon^2f_1\hat U^1f_1^*+\upsilon^1f_2\hat U^2f_2^*\\
    &=(\upsilon^2\tilde 1+f_2\hat U^2f_2^*)(\upsilon^1\tilde 1+f_1\hat U^1f_1^*).
\end{align*}
\end{lemma}

\textbf{Proof:} We have a natural isomorphism $\mathfrak K\mathfrak A\cong \mathfrak A\otimes \mathfrak{K}\mathds F$, and the formula is induced from the action of $\mathscr L\underline 1$ on $\mathfrak{K}\mathds F$ by conjugation. The second statement follows from $f^{*}_1f_2= 0=f^{*}_2f_1$.$\blacksquare$\\

\begin{theorem}
    Let $\mathfrak A\in \texttt C^*\texttt{Alg}_{\mathds F}$. We have a natural $\mathscr L$-permutative category structure on $\texttt{pr}_{\widetilde{\mathfrak{KA}}}$ defined as
    \begin{align*}
    \oplus_{\underline a} \boldsymbol U^x
    &:=\mqty[ \begin{cases}
        \upsilon^x_{\mu\nu}\tilde 1+U^x_{\mu\nu},&x=x'\\
        0,& x\neq x'
    \end{cases}
    ],&
    \tau^{\oplus,\sigma}_{\underline a,\langle \boldsymbol P^x\rangle}
    &:=\boldsymbol{\sigma (\underline p^x)}\oplus_{\underline a}\boldsymbol P^x;\\
    \otimes_f \boldsymbol U^i
    &:= \mqty[ \prod_{\underline m}\upsilon^{i}_{\mu^i\nu^i}\tilde 1
    +\sum_{\underline m} \left(\prod_{i'\neq i} \upsilon^{i'}_{\mu^{i'}\nu^{i'}}\right)
    f_i\hat U^{i}_{\mu^i\nu^i}f_i^*
    ]
    ,&
    \tau^{\otimes,\sigma}_{f,\langle \boldsymbol P^i\rangle}
    &:=\boldsymbol{\sigma \langle\underline p^i\rangle}\otimes_{f\cdot\sigma}\boldsymbol P^i;\\
    \noalign{\centering$\delta_{f,\langle \underline a^i\rangle,\langle \boldsymbol P^{ix}\rangle}:=\boldsymbol{\delta_{\underline m,\langle \underline a^i\rangle,\langle \underline p^{ix}\rangle}}\otimes_f\oplus_{\underline x^i}\boldsymbol P^{ix}.$}
    \end{align*}
\end{theorem}

\textbf{Proof:} The $\mathscr C$-category structure is induced by the direct sum permutative structure. That $\otimes$ is a well-defined functor, which in particular preserve projections and partial isometries, follows from its compatibility with composition and adjunction:
\begin{align*}
    (\otimes_f \boldsymbol U^{i1})(\otimes_f \boldsymbol U^{i2})&=\mqty[\textstyle
    \left(\prod_{\underline m}\upsilon^{i1}_{\mu^i\xi^i}\upsilon^{i2}_{\xi^i\nu^i}\right)\tilde 1
    +\sum_{\underline m}\left(\prod_{i'\neq i}\upsilon^{i'1}_{\mu^{i'}\xi^{i'}}\upsilon^{i'2}_{\xi^{i'}\nu^{i'}}\right)
    f_i(\upsilon^{i1}_{\mu^i\xi^i}\hat U^{i2}_{\xi^i\nu^i}
    +\hat U^{i1}_{\mu^i\xi^i}\upsilon^{i2}_{\xi^i\nu^i}
    +\hat U^{i1}_{\mu^i\xi^i}\hat U^{i2}_{\xi^i\nu^i})f_i^*
    ]\\
    &=\otimes_f \boldsymbol U^{i1}\boldsymbol U^{i2};\\
    (\otimes_f \boldsymbol U^i)^*
    &=\mqty[
    \left(\prod_{\underline m} \bar \upsilon^{i}_{\mu^i\nu^i} \right)\tilde 1
    +\sum_{\underline m}\left(\prod_{i'\neq i} \bar \upsilon^{i'}_{\mu^{i'}\nu^{i'}}\right)
    f_i\hat U^{i,*}_{\mu^i\nu^i}f_i^*
    ]\\
    &=\otimes_f \boldsymbol U^{i,*}.
\end{align*}

The braiding is of the correct form:
    \begin{align*}
        \boldsymbol{\sigma\langle \underline p^{i}\rangle}(\otimes_{f\cdot \sigma}\boldsymbol{P}^{i})\boldsymbol{\sigma\langle \underline p^{i}\rangle}^*
        &=\mqty[
        \left(\prod_{\underline m}\pi^{i}_{\mu^{\sigma i}\nu^{\sigma i}}\right)\tilde 1
        +\sum_{\underline m}
        \left(\prod_{i'\neq i}\pi^{\sigma i'}_{\mu^{\sigma i'}\nu^{\sigma i'}}\right)
        f_{\sigma i}\hat P^{i}_{\mu^{\sigma i}\nu^{\sigma i}}f_{\sigma i}^*
        ]\\
        &=\otimes_f\boldsymbol P^{\sigma^{-1}i}.
\end{align*}

The $\mathscr L$-structure coherence conditions hold:
\begin{align*}
    \otimes_f\otimes_{g^i}\boldsymbol{U}^{ij}
    &=\mqty[
        \left(\prod_{\Sigma_{\underline m}\underline n^i}\upsilon^{ij}_{\mu^{ij}\nu^{ij}}\right)\tilde 1
        +\sum_{\Sigma_{\underline m}\underline n^i}
        \left(\prod_{i'j'\neq ij}\upsilon^{i'j'}_{\mu^{i'j'}\nu^{i'j'}}\right)f_ig^i_j
        \hat U^{ij}_{\mu^{ij}\nu^{ij}}g^{i,*}_jf^{*}_i
    ]\\
    &=\otimes_{f\langle g^i\rangle}\boldsymbol{U}^{ij};\\
    \otimes_{\text{id}}\boldsymbol U
    &=\left[\upsilon_{\mu\nu}\tilde 1
    +\text{id}\hat U_{\mu\nu}\text{id}^*\right]\\
    &=\boldsymbol U;\\
    \tau^{\otimes,\sigma}_{f,\langle \boldsymbol P^{\sigma^{\prime,-1}i}\rangle}
    \tau^{\otimes,\sigma' }_{f\cdot\sigma,\langle \boldsymbol P^{i}\rangle}
    &=\left(\boldsymbol{\sigma\langle\underline p^{\sigma^{\prime-1}i}\rangle}\boldsymbol{\sigma'\langle\underline p^{i}\rangle}\right)\otimes_{f\cdot\sigma\sigma'}\boldsymbol P^{i}\\
    &=\tau^{\otimes,\sigma\sigma'}_{f,\langle \boldsymbol P^{i}\rangle};\\
    \tau^{\otimes,\text{id}}_{f,\langle \boldsymbol P^{i}\rangle}&=\boldsymbol{\text{id}\langle\underline p^i\rangle}\otimes_f \boldsymbol P^i\\
    &=\otimes_f\boldsymbol P^i.
\end{align*}

That $0_{\underline 0}$ is an absorbing element is due to $ M_{\underline 0}\widetilde{\mathfrak{KA}}=\{0_{\underline 0}\}$, and the coherence follows from the following function equations:
\begin{align*}
    \delta_{\underline m,\langle \underline a^{i}\rangle,\langle \underline p^{i(\tau^{i,-1}x)}\rangle}\left(\Pi_{\underline m}\sigma^i( \underline p^{ix})\right)
    &=\left(\left(\Pi_{\underline m}\sigma^i\right)(\underline p^{ix^i})\right)\delta_{\underline m,\langle x^{i}\rangle,\langle \underline p^{ix}\rangle};\\
    \delta_{\underline m,\langle \underline a^{\sigma^{-1}i}\rangle,\langle \underline p^{(\sigma^{-1}i)x}\rangle}\left(\sigma\langle \Sigma_{\underline a^{i}}\underline p^{ix}\rangle\right)
    &=\left(\left(\sigma\langle \underline a^{i}\rangle\right)(\underline p^{(\sigma^{-1}i)x^{\sigma^{-1}i}})\right)\left(\Sigma_{\Pi_{\underline m}\underline a^{i}}\sigma\langle \underline p^{ix^{i}}\rangle\right)\delta_{\underline m,\langle\underline a^{i}\rangle,\langle \underline p^{ix}\rangle};\\
    \left( \delta_{\underline m,\langle  \underline a^i\rangle,\langle \underline b^{ix}\rangle}(\Pi_{\underline m}\underline p^{ix^iy^i})\right)\delta_{\underline m,\langle \Sigma_{\underline a^i}\underline b^{ix}\rangle,\langle \underline p^{ixy}\rangle}
    &=\left(\Sigma_{\Pi_{\underline m}\underline a^i} \delta_{\underline m,\langle \underline b^{ix^i}\rangle,\langle \underline p^{ix^iy}\rangle}\right)\delta_{\underline m,\langle \underline a^i\rangle,\langle \Sigma_{\underline b^{ix}}\underline p^{ixy}\rangle}.\blacksquare
\end{align*}

\section{The special $\hat{\mathscr L}\wr \mathscr F$-spaces $\lVert \mathfrak A\rVert$}

In this section we will give the basic definitions of categories of operators and of ring operators, as well as how to construct examples out of operads and operad pairs. We will then construct special $\hat{\mathscr L}\wr\mathscr F$-spaces out of $\mathscr L$-permutative categories by a method that closely parallels the construction of $\mathscr F\wr\mathscr F$-spaces out of bipermutative categories described in \cite{MayMultLoopSpcTh,Ma09b}. Street's rectification of lax funtors in \cite{street1972two} is central to this construction.

\subsection{Categories of ring operators}

If we have a sequence $\langle\underline m^d_*\rangle\in\Pi_{\underline r_*}\mathscr F$ we will denote a variable in $\underline m^d_*$ as $i_d$. Let $\Pi\subset\mathscr F$ be the subcategory with the same objects and 
$$
    \Pi(\underline m^0_*,\underline m^1_*):=\{\psi\in\mathscr F(\underline m^0_*,\underline m^1_*)\mid\forall i_1\in \underline m^1: \lvert \psi^{-1}i_1\rvert\leq 1\},
$$
ie the functions in $\Pi$ that are injective when restricted to the elements that are not mapped to the base-points. Let also $\Upsilon\subset \Pi$ be the subcategory with the same objects and morphisms 
$$
    \Upsilon(\underline m^0_*,\underline m^1_*):=\{\psi\in\Pi(\underline m^0_*,\underline m^1_*)\mid\forall i_1\in \underline m^1: \lvert \psi^{-1}i_1\rvert= 1\},
$$
ie the functions in $\Upsilon$ are bijective when restricted to the elements that are not mapped to the base-points.

There is a natural inclusion of $\mathbb F^{\text{inj}}$ into $\Pi$, which induces an inclusion of $\mathbb F^{\text{bij}}$ into $\Upsilon$, by adjoining the initial base element 0. 

For all $(\psi^2,\psi^1)\in\mathscr F(\underline m^1_*,\underline m^2_*)\times \mathscr F(\underline m^0_*,\underline m^1_*)$ define for each $i_2\in\underline m^2$ the bijection
$$
    \sigma_{i_2}(\psi^2,\psi^1)\in\mathbb F^{\text{bij}}((\psi^2\psi^1)^{-1}i_2, \Sigma_{\psi^{2,-1} i_2}\psi^{1,-1} i_1),
    \qquad
    \sigma_{i_2}(\psi^2,\psi^1)i_0:=(\psi^1 i_0,i_0).
$$
which is the inverse of the projection $(i_1,i_0)\mapsto i_0$. 

If $i_2\not\in\text{Im }\psi^2\psi^1$ then we set $(\psi^2\psi^1)^{-1}i_2=\Sigma_{\psi^{2,-1} i_2}\psi^{1,-1} i_1:=\underline 1$ and $\sigma_{i_2}(\psi^2,\psi^1):=\text{id}_{\underline 1}$. For all $\underline m_*\in\mathbb F_*$ there is $\phi_{\underline m_*}\in \mathscr F(\underline m_*,\underline 1_*)$ with $\phi_{\underline m_*} i:=1$ if $i>0$. These along with the morphisms in $\Pi$ generate $\mathscr F$. Note that $\wedge_{\underline m} \phi_{\underline p_*^i} = \phi_{\wedge_{\underline m} \underline p_*^i}$. We also have for all $i\in\underline m$ the morphism $\delta_{\underline m_*,i}\in \Upsilon(\underline m_*,\underline 1_*)$ with $\delta_{\underline m_*,i} i':=\delta_{ii'}$, where on the right we are using the Kronecker delta.

\begin{definition}
    A \textit{category of operators} is a topological category $\mathscr D$ with $\mathbb F_*$ as the set of objects and equipped with a factorization 
    $$
        \xymatrix@=0.6cm{\Pi\ar@{^{(}->}[r]^{\eta}&\mathscr D\ar@{->>}[r]^{\epsilon} &\mathscr F}
    $$ of the inclusion $\Pi\hookrightarrow \mathscr F$ with $\eta$ an injection and $\epsilon$ a surjection.
\end{definition}

The morphism spaces $\mathscr D(\underline m^0_*,\underline m^1_*)$ are spaces of operations with $m^0$ inputs and $m^1$ outputs. Operads with non-empty underlying spaces induce categories of operators.

\begin{definition}
    Let $\mathscr G$ be an operad such that $\mathscr G\underline m\neq \emptyset$ for all $\underline m\in \mathbb F$. The category of operators $\hat{\mathscr G}$ has the same objects as $\mathscr F$ and its morphism spaces are
    $$
        \hat{\mathscr G}(\underline m^0_*,\underline m^1_*)
        :=\coprod_{\mathscr F(\underline m^0_*,\underline m^1_*)}\prod_{\underline m^1}\mathscr G(\epsilon\hat\alpha^{-1}i_1).
    $$
    We may denote each particular morphism as $\hat f=(\epsilon\hat f,\langle f^{i_1}\rangle)$, so that with this notation the compositions, units and structural maps are
    \begin{gather*}
        \hat f^2\hat f^1:=(\epsilon\hat  f^2\epsilon\hat f^1,\langle  f^{2i_2} \langle  f^{1i_1}\rangle\cdot \sigma_{i_2}(\epsilon\hat  f^2,\epsilon\hat f^1)\rangle),
        \qquad
        \text{id}_{\underline m_*}:=(\text{id}_{\underline m_*},\langle\textbf{id}\rangle);\\
        \eta (\psi):=\left(\psi,\left\langle\begin{cases}
        \ast,&\lvert \psi^{-1}i_1\rvert=0\\
        \textbf{id},&\lvert \psi^{-1}i_1\rvert=1
        \end{cases}\right\rangle\right),
        \qquad
        \epsilon(\hat f):=\epsilon\hat  f.
    \end{gather*}
\end{definition}

In $\hat{\mathscr L}$ for instance the composition is
$$
    \hat f^2\hat f^1
    =(\epsilon\hat f^2 \epsilon \hat f^1,\langle\langle f^{2i_2}_{\epsilon\hat f^1 i_0}f^{1(\epsilon\hat f^1i_0)}_{i_0}\rangle\rangle).
$$

We have $\hat{\mathscr C }\cong\mathscr F$. For $\mathds 1$ the reduced operad with $\mathds 1\underline 1=\{\textbf{id}\}$ and $\mathds 1\underline n=\emptyset$ for $n>1$ we have $\hat{\mathds 1}=\Pi$. For $\mathds 1^\circ$ the non-reduced operad with $\mathds 1^\circ\underline 1=\{\textbf{id}\}$ and $\mathds 1^\circ\underline n=\emptyset$ for $n\neq 1$ we have $\hat{\mathds 1}^\circ=\Upsilon$.

As with operads, distributivity laws between categories of operators can be encoded in actions.

\begin{definition}
     Let $\mathscr D$ and $\mathscr K$ be categories of operators. An \textit{action} $\ltimes$ of $\mathscr K$ on $\mathscr D$ consists of functors $\hat f\ltimes:\mathscr D^{\underline m^0}\rightarrow \mathscr D^{\underline m^1}$ for all $\hat f\in\mathscr K(\underline m^0_*,\underline m^1_*)$ such that:
     \begin{enumerate}[(i)]
        \item If $\langle\underline a_*^{i_0}\rangle\in\mathscr D^{\underline m^0}$ then $\hat f\ltimes\langle\underline a_*^{i_0}\rangle=\langle \wedge_{\epsilon\hat f^{-1}i_1}\underline a_*^{i_0}\rangle$;
         
        \item If $\langle \iota^{i_0}\rangle\in \Pi^{\underline m^0}(\langle\underline a^{i_00}_*\rangle,\langle\underline a^{i_01}_*\rangle)$ then 
        $\hat f\ltimes\langle \eta \iota^{i_0}\rangle=\langle \eta_{\mathscr D}(\wedge_{\epsilon \hat f^{-1}i_1}\iota^{i_0})\rangle$;
        
        \item If $\langle \hat \alpha^{i_0}\rangle\in \mathscr D^{\underline m^0}(\langle\underline a^{i_00}_*\rangle,\langle\underline a^{i_01}_*\rangle)$ then 
        $\epsilon ^{\underline m^1}(\hat f\ltimes\langle \hat\alpha^{i_0}\rangle)=\langle \eta (\wedge_{\epsilon  \hat f^{-1}i_1}\epsilon \hat \alpha^{i_0})\rangle$;
         
         \item If $\psi\in\Pi(\underline m^0_*, \underline m^1_*)$ then 
         $\eta \psi\ltimes \langle \hat \alpha^{i_0}\rangle:=\langle \hat \alpha^{\psi^{-1}i_1}\rangle$;
         
         \item If $(\hat f^2,\hat f^1)\in\mathscr K(\underline m^1_*,\underline m^2_*)\times \mathscr K(\underline m^0_*,\underline m^1_*)$ then the isomorphism $\langle\eta \sigma_{i_2}(\epsilon \hat f^2,\epsilon \hat f^1)\rangle$ specifies a natural isomorphism $(\hat f^2\hat f^1)\ltimes \Rightarrow (\hat f^2\ltimes)( \hat f^1\ltimes)$.
     \end{enumerate}
     In (i) if $(\epsilon \hat f)^{-1}i_1=\emptyset$ then the $i_1$-th coordinate is $\underline 1_*$, and in (ii)-(iv) this implies the $i_1$-th coordinate is $\text{id}_{\underline 1_*}$. 
     
     A pair of categories of operators equipped with an action $(\mathscr D,\mathscr K,\ltimes)$ will be refered to as a \textit{categories of operators pair}.
\end{definition}

In particular $\Pi$ acts on any category of operators $\mathscr D$ by the formulas in (i) and (iv), and a category of operators $\mathscr K$ acts on both $\Pi$ and $\mathscr F$ by the formulas in (i), (ii) and (iii). 

We have natural embeddings 
$$\iota:\mathscr P \underline m\rightarrow \hat{\mathscr P}(\underline m_*,\underline 1_*),
\qquad 
\iota\alpha:=(\phi_{\underline m_*},\langle\alpha\rangle),
$$
and operadic composition of $\mathscr P$ is recovered from categorial composition in $\hat{\mathscr P}$. An action of $\mathscr G$ on $\mathscr P$ determines and is determined by an action of $\hat{\mathscr G}$ on $\hat{\mathscr P}$. From an action of $\hat{\mathscr G}$ on $\hat{\mathscr P}$ the formula 
$\iota(f\ltimes\langle \alpha^i\rangle):=(\iota f)\ltimes\langle \iota\alpha^i\rangle$ determines an action of $\mathscr G$ on $\mathscr P$. Conversely an action of $\mathscr G$ on $\mathscr P$ induces an action of $\hat{\mathscr G}$ on $\hat{\mathscr P}$ by the formula
$$
\hat f\ltimes\langle \hat \alpha^{i_0}\rangle=\langle (\wedge_{\epsilon\hat f^{-1}i_1}\epsilon \hat\alpha^{i_0},\langle f^{i_1}\ltimes \langle \alpha^{i_0 a_{i_01}}\rangle\rangle) \rangle.
$$

Wreath products allow us to condense the information contained in a pair of categories of operators equipped with an action into a single category.

\begin{definition}
    Let $(\mathscr D,\mathscr K,\ltimes)$ be a categories of operator pair. The \textit{wreath product} $\mathscr K\wr \mathscr D$ is the category with $\Sigma_{\mathbb F_*}\mathbb F_*^{\underline m}$ as objects, i.e. pairs of the form $(\underline m_*,\langle \underline a_*^i\rangle)$. There is a unique object $(\underline 0_*,\ast)$ when $m=0$. The morphism spaces are 
$$
    \mathscr K\wr \mathscr D((\underline m^0_*,\langle \underline a_*^{0i_0}\rangle),(\underline m^1_*,\langle \underline a_*^{1i_1}\rangle))
    :=\coprod_{\mathscr F(\underline m^0_*,\underline m^1_*)}\epsilon_{\mathscr K}^{-1}\psi\times\prod_{\underline m^1}\mathscr D(\wedge_{\psi^{-1}i_1}\underline a_*^{0i_0},\underline a_*^{1i_1}).
$$
The compositions and identities are 
\begin{align*}
    (\hat f^2,\langle \hat \alpha^{2i_2}\rangle)(\hat f^1,\langle \hat \alpha^{1i_1}\rangle)
    &:=(\hat f^2\hat f^1,\langle \hat \alpha^{2i_2}\rangle( \hat f^2\ltimes\langle \hat \alpha^{1i_1}\rangle)\langle\eta (\sigma_{i_2}(\epsilon \hat f^2,\epsilon \hat f^1)\langle \underline a^{0i_0}\rangle)\rangle),\\
    \text{id}_{(\underline m_*,\langle \underline a_*^i\rangle)}
    &:=(\text{id}_{\underline m_*},\langle \text{id}_{\underline a^i_*}\rangle).
\end{align*}
\end{definition}

In particular in $\hat{\mathscr L}\wr \mathscr F$ we have
\begin{gather*}
    (\hat f^2,\langle \kappa^{2i_2}\rangle)(\hat f^1,\langle  \kappa^{1i_1}\rangle)=
    ( \hat f^2 \hat f^1,\langle \kappa^{2i_2} (\wedge_{\epsilon\hat f^{2,-1}i_2}\kappa^{1i_1})(\sigma_{i_2}(\epsilon \hat f^2,\epsilon \hat f^1)\langle\underline a^{0i_0}\rangle)\rangle).
\end{gather*}

With the wreath product we can define categories of ring operators in a similar way that categories of operators are defined.

\begin{definition}
    A \textit{category of ring operators} is a topological category $\mathscr J$ with $\Sigma_{\mathbb F_*}\mathbb F^{\underline m}_*$ as the set of objects and equipped with a factorization 
    $$
        \xymatrix@=0.6cm{\Pi\wr\Pi\ar@{^{(}->}[r]^(0.55){\eta}&\mathscr J\ar@{->>}[r]^(0.4){\epsilon} &\mathscr F\wr\mathscr F}
    $$
    of the inclusion $\Pi\wr\Pi\hookrightarrow \mathscr F\wr\mathscr F$ with $\eta$ an injection and $\epsilon$ a surjection.
\end{definition}

For every operad pair $(\mathscr P,\mathscr G,\ltimes)$ the wreath product $\hat{\mathscr G}\wr\hat{\mathscr P}$ is a category of ring operators.

Functors from a category of ring operators to the category of topological spaces satisfying some special conditions behave like semirings up to homotopy.
\begin{definition}
    Let $\mathscr J$ be a category of ring operators. A \textit{$\mathscr J$-space} $Z$ is a continuous functor $Z:\mathscr J\rightarrow\texttt{Top}$. It is \textit{semispecial} if
    \begin{enumerate}
        \item $Z(\underline 0_*,\ast)$ is aspherical.
        
        \item The maps $\langle Z(\eta\delta_{\underline m_*,i},\text{id}_{\underline a^i_*}) \rangle:Z(\underline m_*,\langle \underline a_*^i\rangle)\rightarrow\Pi_{\underline m}Z(\underline 1_*, \underline a_*^i)$ are equivalences.
    \end{enumerate}
        If $Z$ is semispecial then $\pi_0 Z(\underline 1_*, \underline 1_* )$ equipped with the product defined as the composition
        $$\xymatrix@=1.8cm{
            \pi_0 Z(\underline 1_*, \underline 1_* )\times\pi_0 Z(\underline 1_*, \underline 1_* )\ar[r]^(0.575){\langle Z(\eta\delta_{\underline 2_*,i},\text{id}_{\underline 1_*})_*\rangle^{-1}}
            &\pi_0 Z(\underline 2_*,\langle \underline 1_*,\underline 1_*\rangle)
            \ar[r]^(0.55){Z(\eta\phi_{\underline 2_*},\text{id}_{\underline 1_*})_*}
            &\pi_0 Z(\underline 1_*, \underline 1_* )
        }$$
        is a monoid. We say $Z$ is \textit{special} if it is semispecial and
    \begin{enumerate}
        \item[3.] $Z(\underline 1_*, \underline 0_*)$ is aspherical.
        
        \item[4.] The maps $\langle Z(\text{id}_{\underline 1_*},\eta\delta_{\underline a_*,x})\rangle: Z(\underline 1_*, \underline a_* )\rightarrow \Pi_{\underline a}Z(\underline 1_*, \underline 1_* )$ are equivalences.
    \end{enumerate}
    If $Z$ is special then $\pi_0 Z(\underline 1_*, \underline 1_* )$ equipped with the sum defined as the composition
    $$\xymatrix@=1.8cm{
            \pi_0 Z(\underline 1_*, \underline 1_* )\times\pi_0 Z(\underline 1_*, \underline 1_* )
            \ar[r]^(0.65){\langle Z(\text{id}_{\underline 1_*},\eta\delta_{\underline 2_*,x})\rangle_*^{-1}}
            &\pi_0 Z(\underline 1_*, \underline 2_* )
            \ar[r]^(0.5){Z(\text{id}_{\underline 1_*},\eta \phi_{\underline 2_*})_*}
            &\pi_0 Z(\underline 1_*, \underline 1_*)
        }$$
    is a semiring. The category of special $\mathscr J$-spaces is denoted as $\mathscr J[\texttt{Top}]$.
\end{definition}

\subsection{Street's rectification of lax functors}

For all $\mathfrak A\in \texttt C^*\texttt{Alg}_{\mathds F}$ we functorialy construct a special $\hat{\mathscr L}\wr \mathscr F$-space out of the $\mathscr L$-permutative category $ \texttt{pr}_{\widetilde{\mathfrak{KA}}}$ by a procedure that closely parallels the construction of special $\mathscr F\wr \mathscr F$-spaces out of bipermutative categories. We construct a lax functor 
$$
    \mathcal F\texttt{pr}_{\widetilde{\mathfrak{KA}}}:\hat{\mathscr L}\wr \mathscr F \leadsto\texttt{Cat}
$$
out of $ \texttt{pr}_{\widetilde{\mathfrak{KA}}}$, rectify it into a functor
$$
\overline{\mathcal F\texttt{pr}_{\widetilde{\mathfrak{KA}}}}:\hat{\mathscr L}\wr \mathscr F \rightarrow\texttt{Cat}
$$
using Street's first construction \cite{street1972two}, and then take geometric realization to get a special $\hat{\mathscr L}\wr \mathscr F$-space
$$\lVert\mathfrak A\rVert:=\lvert\overline{\mathcal F\texttt{pr}_{\widetilde{\mathfrak{KA}}}}\rvert:\hat{\mathscr L}\wr \mathscr F \rightarrow\texttt{Top}.
$$
This method applies to any $\mathscr L$-permutative category, but we focus on the $ \texttt{pr}_{\widetilde{\mathfrak{KA}}}$ since they are our main examples of interest.

In higher category theory a \textit{lax functor} is a morphism between bicategories that preserves composition and identities of 1-morphisms only up to coherently specified 2-isomorphisms. We will not need the full theory of bicategories, only the following characterization of lax functors with domain a category (considered as a bicategory with only identity 2-cells) and codomain $\texttt{Cat}$, the bicategory of small categories, functors and natural transformations.

\begin{definition}
    For $\mathcal C$ a category a \textit{lax functor} $F:\mathcal C\leadsto\texttt{Cat}$ consists of:
    \begin{itemize}
        \item For each $X\in\mathcal C$ a category
        $$FX\in\texttt{Cat};$$
        
        \item For each $\psi\in \mathcal C(X^0,X^1)$ a functor 
        $$F\psi\in\texttt{Cat}(FX^0, FX^1);
        $$
        
        \item For each $(\psi^2,\psi^1)\in \mathcal C(X^1,X^2)\times\mathcal C(X^0,X^1)$ a natural transformation
        $$\xymatrix@=2cm{
            FX^0 \rtwocell^{F\psi^2F\psi^1}_{F(\psi^2\psi^1)}{\quad\ \ \omega^{\psi^2,\psi^1}} & FX^2
        };$$
        
        \item For each $X\in\mathcal C$ a natural transformation
        $$\xymatrix@=1.25cm{
           FX\rtwocell^{\text{id}_{FX}}_{F\text{id}X}{\ \ \omega^X}&FX
        };$$
    \end{itemize}
    satisfying the coherence condition that the following diagrams commute:
    \begin{gather*}\xymatrix{
        F\psi^3F\psi^2F\psi^1
        \ar@{=>}[r]^(0.525){\text{id}_{F\psi^3}\bullet \omega^{\psi^2,\psi^1}}
        \ar@{=>}[d]_{\omega^{\psi^3,\psi^2}\bullet \text{id}_{F\psi^1}}
        &F\psi^3F(\psi^2 \psi^1)
        \ar@{=>}[d]^{\omega^{\psi^3,\psi^2\psi^1}}\\
        F(\psi^3\psi^2) F\psi^1
        \ar@{=>}[r]_{\omega^{\psi^3\psi^2,\psi^1}}
        &F(\psi^3\psi^2\psi^1)
    },\quad
    \xymatrix@C=0.95cm{
        &F \psi
        \ar@{=>}[dl]_{\omega^{X^1} \bullet\text{id}_{F\psi}}
        \ar@{=}[d]
        \ar@{=>}[dr]^{\text{id}_{F\psi}\bullet\omega^{X^0}}&\\
        F\text{id}_{X^1} F\psi
        \ar@{=>}[r]_(0.6){\omega^{\text{id}_{X^1},\psi}}
        &F\psi
        & F\psi F\text{id}_{X^0}
        \ar@{=>}[l]^(0.6){\omega^{\psi,\text{id}_{X^0}}}
    }.\end{gather*}
\end{definition}

\begin{theorem}
    For each $\mathcal C\in (\mathscr C,\mathscr L)[\texttt{Cat}]$ we can construct the lax functor
\begin{gather*}
    \mathcal{FC}: \hat{\mathscr L}\wr\mathscr F\leadsto\texttt{Cat},\\
    \mathcal{FC}(\underline m_*,\langle \underline a^i_*\rangle)
    := \mathcal C^{\Sigma_{\underline m}\underline a^i},
    \qquad\mathcal{FC}(\hat f,\langle \kappa^{i_1}\rangle)
    \langle P^{i_0x_0}\rangle
    :=\langle\oplus_{\kappa^{i_1,-1}x_1}\otimes_{f^{i_1}}P^{i_0x_0^{i_0}}\rangle,\\
    \xymatrix@=1cm{
    \langle\oplus_{\kappa^{2i_2,-1}x_2}\otimes_{f^{i_2}} \oplus_{\kappa^{1i_1,-1}x_1^{i_1}}\otimes_{f^{i_1}} P^{i_0x_0^{i_0}}\rangle
    \ar@{=>}@/_1cm/@<-2cm>[ddd]_{\omega^{(\hat f^2,\langle \kappa^{2i_2}\rangle),(\hat f^1,\langle \kappa^{1i_1}\rangle)}_{\langle P^{i_0x_0}\rangle}}
    \ar@{=>}[d]^(0.5){\left\langle\oplus_{\kappa^{2i_2,-1}x_2}\delta_{f^{i_2},\langle \kappa^{1i_1,-1}x_1^{i_1}\rangle,\langle\otimes_{f^{i_1}}P^{i_0x_0^{i_0}}\rangle}\right\rangle}
    \\
    \langle\oplus_{(\kappa^{2i_2}(\wedge_{\epsilon\hat f^{2,-1}i_2}\kappa^{1i_1}))^{-1}x_2}\otimes_{f^{i_2}\langle f^{i_1}\rangle} P^{i_0x_0^{i_0}} \rangle
    \ar@{=>}[d]^{\left\langle\oplus_{(\kappa^{2i_2}(\wedge_{\epsilon\hat f^{2,-1}i_2}\kappa^{1i_1}))^{-1}x_2}\ \tau^{\otimes,\sigma_{i_2}(\epsilon\hat f^2,\epsilon\hat f^1)^{-1}}_{\langle f^{i_2}_{\epsilon\hat f i_0} f^{\epsilon\hat f i_0}_{i_0}\rangle,\langle P^{i_0x_0^{i_0}}\rangle}\right\rangle}\\
    \langle\oplus_{(\kappa^{2i_2}(\wedge_{\epsilon\hat f^{2,-1}i_2}\kappa^{1i_1}))^{-1}x_2}\otimes_{\langle f^{i_2}_{\epsilon\hat fi_0} f^{\epsilon\hat fi_0}_{i_0}\rangle}P^{i_0x_0^{i_0}} \rangle
    \ar@{=>}[d]^(0.55){\left\langle\tau^{\oplus,\sigma_{i_2}(\epsilon\hat f^2,\epsilon\hat f^1)^{-1}\langle\underline a^{0i_0}\rangle
    }_{(\kappa^{2i_2}(\wedge_{\epsilon\hat f^{2,-1}i_2}\kappa^{1i_1})(\sigma_{i_2}(\epsilon\hat f^2,\epsilon\hat f^1)\langle\underline a^{0i_0}\rangle))^{-1}x_2,\langle\otimes_{ \langle f^{i_2}_{\epsilon\hat fi_0} f^{\epsilon\hat fi_0}_{i_0}\rangle} P^{i_0x_0^{i_0}} \rangle} \right\rangle}\\
    \langle\oplus_{(\kappa^{2i_2}(\wedge_{\epsilon\hat f^{2,-1}i_2}\kappa^{1i_1})(\sigma_{i_2}(\epsilon\hat f^2,\epsilon\hat f^1)\langle \underline a^{0i_0}\rangle))^{-1}x_2}\otimes_{\langle f^{i_2}_{\epsilon\hat f i_0} f^{\epsilon\hat fi_0}_{i_0}\rangle} P^{i_0x_0^{i_0}} \rangle},\quad
    \xymatrix@=2.5cm{\langle  P^{ix}\rangle
    \ar@{=>}[r]^{\omega^{(\underline m_*,\langle\underline a_*^i\rangle)}_{\langle P^{ix}\rangle}:=\text{id}_{\langle  P^{ix}\rangle}=\langle P^{ix}\rangle}
    &\langle P^{ix}\rangle}.
\end{gather*}
\end{theorem}

\textbf{Proof:} That the lax functor cohesion diagrams commute follows from the cohesion of the $\mathscr L$-permutative category structure.$\blacksquare$

We are particularly interested in the lax functors $\mathcal F\texttt{pr}_{\widetilde{\mathfrak{KA}}}$ for $\mathfrak A\in \texttt C^*\texttt{Alg}_{\mathds F}$.

\begin{definition}
    Given a lax functor $F:\mathcal C\leadsto \texttt{Cat}$ \textit{Street's first construction} yields the functor
\begin{align*}
    \noalign{\centering$\overline F:\mathcal C\rightarrow\texttt{Cat}$}
    \overline FX&:=
    \begin{cases}
    \text{Ob}\ \overline FX:=\coprod_{\mathcal C_{/X}} \text{Ob}\ FX^1;\\
    \overline FX((\phi^0,P^1),(\phi^1,P^2))
    :=\coprod_{\mathcal C_{/_X}(\phi^1,\phi^0)}F X^1(P^1, F\psi P^2);\\
    (\psi^2,U^2)(\psi^1,U^1):=(\psi^1\psi^2,\omega^{\psi^1,\psi^2}_{P^2}(F\psi^1 U^2)U^1);\\
    \text{id}_{(\phi,P)}:=(\text{id}_{X^1},\omega^{X^1}_P).
    \end{cases},\quad
    \overline F(\chi:X\rightarrow Y)
    &:=\begin{cases}
        \overline F\chi (\phi,P):=(\chi\phi,P);\\
        \overline F\chi (\psi,U):=(\psi,U).
    \end{cases}
\end{align*}
\end{definition}

\begin{figure}
    \centering
    $$\xymatrix{
&X&\\
X^{1}\ar[ur]^{\phi^0}
&X^{2}\ar[u]|{\phi^1}
\ar[l]^{\psi^1}
&X^{3}\ar[ul]_{\phi^2}
\ar[l]^{\psi^2}\\
P^1\ar[d]_{U^1}
&P^2 \ar@{|->}[dl]
\ar[d]_{U^2}
&P^3 \ar@{|->}[dl]\ar@/^0.5cm/@{|->}[dddll]\\
F\psi^1 P^2\ar[d]_{F\psi^1U^2}
&F\psi^2 P^3\ar@{|->}[dl]&\\
F\psi^1F\psi^2 P^3\ar[d]_{\omega_{P^3}^{\psi^1,\psi^2}}&&\\
F(\psi^1\psi^2)P^3&&
}\qquad
\xymatrix{
&X&\\
X^1\ar[ur]^{\phi}&&X^1\ar[ll]^{\text{id}_{X^1}}\ar[ul]_{\phi}\\
P\ar[d]_{\omega^{X^1}_P}&&P\ar@{|->}[dll]\\
F\text{id}_{X^1}P&&
}
$$
    \caption{Diagrammatic representation of compositions and identities in $\overline F X$}
    \label{fig:com and id in Ovln FC}
\end{figure}
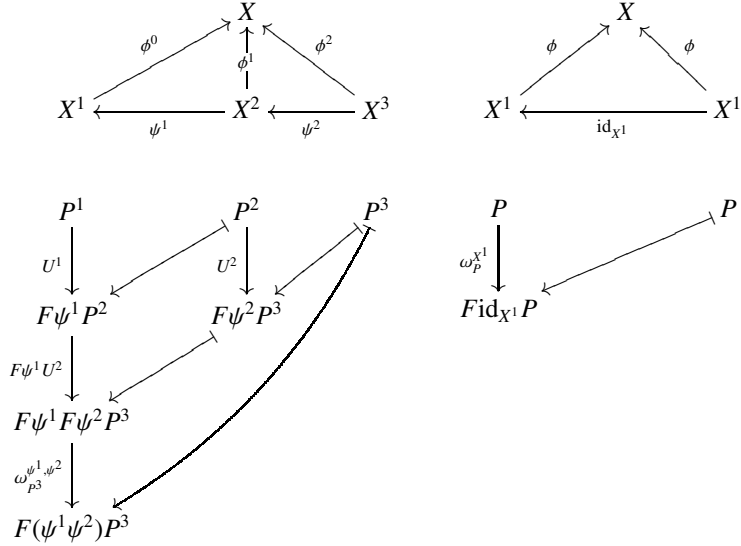

A morphism from $(\phi^0,P^1)$ to $(\phi^1,P^2)$ involves a morphism in $\mathcal C_{/X}$ from $\phi^1$ to $\phi^0$, an inversion of the orientation one might naively expect. Applying this construction on the lax functor $\mathcal F\texttt{pr}_{\widetilde{\mathfrak{KA}}}$ for $\mathfrak A\in \texttt{C$^*$Alg}$ we get
\begin{align*}
    \text{Ob}\ \overline{\mathcal F\texttt{pr}_{\widetilde{\mathfrak{KA}}}}(\underline m^0,\langle \underline a^{0i_0}\rangle)
    &=\coprod_{\hat{\mathscr L}\wr \mathscr F_{/(\underline m^0,\langle \underline a^{0i_0}\rangle)}}(\mathcal P_\infty\widetilde{\mathfrak{KA}})^{\Sigma_{\underline{m}^1}\underline a^{1i_1}},\\
    \overline{\mathcal F\texttt{pr}_{\widetilde{\mathfrak{KA}}}}(\underline m^0_*,\langle \underline a^{0i_0}_*\rangle)(((\hat f^0,\langle \kappa^{0i_0}\rangle), \langle\boldsymbol P^{i_1x_1}\rangle),((\hat f^1,\langle \kappa^{1i_1}\rangle), \langle\boldsymbol P^{i_2x_2}\rangle))
    &=\hspace{-1.5cm}
    \coprod_{\hspace{1.25cm}
    \hat{\mathscr L}\wr \mathscr F_{/(\underline m^0_*,\langle \underline a^{0i_0}_*\rangle)}
    ((\hat f^1,\langle \kappa^{1i_1}\rangle)
    ,(\hat f^0,\langle \kappa^{0i_0}\rangle))}
    \hspace{-1.35cm}
    \texttt{pr}_{\widetilde{\mathfrak{KA}}}^{\Sigma_{\underline m^1}\underline a^{1i_1}}
    (\langle P^{i_1x_1}\rangle,\langle\oplus_{\lambda^{i_1,-1}x_1}\otimes_{g^{i_1}}\boldsymbol P^{i_2x_2^{i_2}}\rangle).
\end{align*}

\begin{lemma}
For all $\mathcal C\in (\mathscr C,\mathscr L)[\texttt{Cat}]$ the $\hat{\mathscr L}\wr\mathscr F$-space $\lvert \overline{\mathcal{FC}}\rvert$ is special.
\end{lemma}

\textbf{Proof:} This follows from the fact that for each $(\underline m^0_*,\langle \underline a^{0i_0}_*\rangle)\in \hat{\mathscr L}\wr\mathscr F$ we have a natural adjunction
\begin{gather*}
    (F\dashv U):\mathcal{FC}(\underline m^0_*,\langle \underline a^{0i_0}_*\rangle)\leftrightharpoons \overline{\mathcal{FC}}(\underline m^0_*,\langle \underline a^{0i_0}_*\rangle);\\
    F\langle P^{i_0x_0}\rangle=(\text{id}_{(\underline m^0_*,\langle \underline a^{0i_0}_*\rangle)}, \langle P^{i_0x_0}\rangle)
    ,\qquad
    U((\hat f,\langle \kappa^{i_0}\rangle),\langle P^{i_1x_1}\rangle)=\langle\oplus_{\kappa^{i_0,-1}x_0}\otimes_{f^{i_0}} P^{i_1x_1^{i_1}}\rangle,
\end{gather*}
with $F$ part of a natural lax transformation $F:\mathcal{FC} \Rightarrow \overline{\mathcal{FC}}$. The details are similar to the ones for the construction on bipermutative categories in \cite[Section 3]{MayMultLoopSpcTh}. $\blacksquare$\\

Putting all the constructions in this section together, and focusing on the Murray-von Neuman $\mathscr L$-permutative categories $\texttt{pr}_{\widetilde{\mathfrak{KA}}}$ we get the following definition.

\begin{definition}
The \textit{classifying Murray-von Neumann special $\hat{\mathscr L}\wr\mathscr F$-space} functor is
$$
\lVert - \rVert :
\texttt{C$^*$Alg}_{\mathds F}
\rightarrow 
\hat{\mathscr L}\wr\mathscr F[\texttt{Top}],
\qquad
\lVert\mathfrak A\rVert:=\lvert\overline{\mathcal F\texttt{pr}_{\widetilde{\mathfrak{KA}}}}\rvert.
$$
\end{definition}

Points of the space $\lVert \mathfrak A\rVert(\underline m_*,\langle\underline a_*^{i_0}\rangle)$ can be diagrammatically denoted as
\begin{gather*}
    [((\hat f^0,\langle\kappa^{0i_0}\rangle),\langle ((\hat f^d,\langle \kappa^{di_d}\rangle),\langle\boldsymbol U^{di_dx_d}\rangle)\rangle,\langle\boldsymbol P^{i_{s+1}x_{s+1}}\rangle),\langle t^d\rangle]=\\
    \left[\xymatrix@=1cm{
    (\underline m^0_*,\langle\underline a_*^{0i_0}\rangle)&
    (\underline m^1_*,\langle\underline a_*^{1i_1}\rangle)
    \ar[l]_(0.5){(\hat f^0,\langle\kappa^{0i_0}\rangle)}&
    &(\underline m^{s+1}_*,\langle\underline a^{(s+1)i_{s+1}}_*\rangle)
    \ar[ll]^(0.525){\langle t^d\rangle}_(0.525){\langle((\hat f^d,\langle\kappa^{di_d}\rangle),\langle\boldsymbol U^{di_dx_d}\rangle)\rangle}
    |(0.525){\ \cdots\ }}
    ,\langle\boldsymbol P^{i_{s+1}x_{s+1}}\rangle\right].
\end{gather*}

Setting
$$
\lVert \mathcal U_{\mathds F}\rVert:=
\lvert \overline{\mathcal{FU}_{\mathds F}}\rvert
$$
the natural inclusion  $\lambda\in \texttt C^*\texttt{Alg}(\mathds C,\widetilde{\mathfrak{KA}})$ induces a natural inclusion 
$$\lambda_*\in \hat{\mathscr L}\wr \mathscr F[\texttt{Top}](\lVert \mathcal U_{\mathds F}\rVert,\lVert \mathfrak A\rVert).$$

\section{The {$E_\infty$}-ring spaces {$\Lambda\lVert\mathfrak A\rVert$}}

It would have been a nice surprise if $\lVert \mathfrak A\rVert(\underline 1_*,\underline 1_*)$ was an $E_\infty$-ring space. Alas, the gods are not that kind.  We must construct spaces of resolutions using the two sided bar construction twice. As described in \cite{MayMultLoopSpcTh,Ma09b} the construction of $(\mathscr H_\infty,\mathscr L)$-spaces out of $\hat{\mathscr L}\wr\hat{\mathscr H}_\infty$-spaces goes through the intermediary category of $(\hat{\mathscr H}_\infty,\hat{\mathscr L})$-spaces.

\subsection{The {$(\hat{\mathscr H}_\infty,\hat{\mathscr L})$}-spaces {$\Lambda''\lVert \mathfrak A\rVert$}}

Pulling back over the functor $\text{id}\wr\epsilon:\hat{\mathscr L}\wr\hat{\mathscr H}_\infty\rightarrow \hat{\mathscr L}\wr \mathscr F$ we have an inclusion of $\hat{\mathscr L}\wr\mathscr F[\texttt{Top}]$ in $\hat{\mathscr L}\wr\hat{\mathscr H}_\infty[\texttt{Top}]$, so in particular $\lVert\mathfrak A\rVert\in \hat{\mathscr L}\wr\hat{\mathscr H}_\infty[\texttt{Top}]$. Since $\mathscr H_\infty$ is a $\Sigma$-cofibrant $E_\infty$-operad this inclusion induces an equivalence of homotopy categories. We will require a monadic description of $\hat{\mathscr L}\wr \hat{\mathscr H}_\infty$-spaces \cite[Construction 5.7]{MayMultLoopSpcTh}. 

For $\mathscr P$ an operad define the monad $(\hat P;\eta,\mu)$ on $\texttt{Top}^\Upsilon$ as
\begin{gather*}
    (\hat P Z')\underline a^0_*:=\int^\Upsilon \hat{\mathscr H}_\infty(\underline a^1_*,\underline a^0_*)\times Z'\underline a^1;\qquad
    \eta_{\underline a^0_*}z:=[\text{id}_{\underline a^0_*},z],
    \quad \mu_{\underline a^0_*}[\hat \alpha^0,[\hat\alpha^1,z]]:=[\hat\alpha^0\hat\alpha^1,z].
\end{gather*}

For an operad pair $(\mathscr P,\mathscr G,\ltimes)$ define the monad $(\hat G\wr\hat P;\eta,\mu)$ on $\texttt{Top}^{\Upsilon\wr\Upsilon}$ as
\begin{gather*}
    \displaystyle ((\hat G\wr\hat P)Z'')(\underline m_*^0,\langle\underline a^{0i_0}_*\rangle):=\int^{\Upsilon \wr \Upsilon}\hat{\mathscr G}\wr\hat{\mathscr P}((\underline m^1_*,\langle\underline a^{1i_1}_*\rangle),(\underline m^0,\langle \underline a^{0i_0}\rangle))\times Z''(\underline m^1,\langle \underline a^{1i_1}\rangle);\\
    \eta_{(\underline m^0_*,\langle\underline a^{0i_0}_*\rangle)}z
    :=[\text{id}_{(\underline m^0_*,\langle\underline a^{0i_0}_*\rangle)},z]
    ,\qquad
    \mu_{(\underline m^0_*,\langle\underline a^{0i_0}_*\rangle)}[(\hat f^0,\langle \hat \alpha^{0i_0}\rangle),[(\hat f^1,\langle \hat \alpha^{1i_1}\rangle),z]]
    :=[(\hat f^0,\langle \hat \alpha^{0i_0}\rangle)(\hat f^1,\langle \hat \alpha^{1i_1}\rangle),z].
\end{gather*}

The restriction functor
$$
    -\!\!\restriction_{\Upsilon \wr \Upsilon}:\hat{\mathscr G}\wr\hat{\mathscr P}[\texttt{Top}]\rightarrow \hat G\wr\hat P[\texttt{Top}^{\Upsilon\wr \Upsilon}]
$$
is an isomorphism. In particular
$$\lVert \mathfrak A\rVert\!\!\restriction_{\Upsilon \wr \Upsilon} \ \in \hat L\wr\hat H_\infty[\texttt{Top}^{\Upsilon \wr \Upsilon}].
$$

We have an adjunction 
\begin{gather*}
    \xymatrix{\texttt{Top}^{\Upsilon}
    \ar@<0.15cm>[r]^{R''}\ar@{}[r]|\top
    &\texttt{Top}^{\Upsilon\wr\Upsilon}
    \ar@<0.15cm>[l]^{L''}};\quad 
    \begin{cases}\textstyle
        R'' Z'(\underline m_*,\langle \underline a^i_*\rangle)
        := \prod_{\underline m} Z'\underline a^i_*\\
        R'' Z'(\psi,\langle \kappa^{i_1}\rangle):= \langle (Z'\kappa^{i_1}) \psi\cdot- \rangle
    \end{cases}\hspace{-0.3cm},
    \ \begin{cases}
        L''Z''\underline a_*:=Z''(\underline 1_*,\underline a_*)\\
        L''Z''\kappa:=Z''(\text{id}_{\underline 1_*},\kappa)
    \end{cases}\hspace{-0.35cm}.
\end{gather*}

For $\mathscr G$ an operad the monad $\hat G\wr\Upsilon$ induces a monad $\tilde G:=L''(\hat G\wr\Upsilon)R''$ on $\texttt{Top}^\Upsilon$. If $(\mathscr P,\mathscr G,\ltimes)$ is an operad pair then the monad $\hat P$ induces a monad, also denoted $\hat P$, on $\tilde G[\texttt{Top}^\Upsilon]$.

\begin{definition}
    Let $(\mathscr P, \mathscr G,\ltimes)$ be an operad pair. A \textit{$(\hat{\mathscr P},\hat{\mathscr G})$-space} is a $\Upsilon$-space $Z'\in\texttt{Top}^\Upsilon$ equipped with a $(\hat G\wr\hat P)$-structure on $R''Z'$. Equivalently it is a $\hat P$-algebra in $\tilde G[\texttt{Top}^{\Upsilon}]$. It is \textit{special} if 
    \begin{enumerate}
        \item $Z'\underline 0_*$ is aspherical;
        
        \item The maps $\langle Z'\delta_{\underline a_*,a}\rangle:Z'\underline a_*\rightarrow \Pi_{\underline a} Z'\underline 1_*$ are equivalences.
    \end{enumerate}
    The category of special $(\hat{\mathscr P},\hat{\mathscr G})$-spaces is denoted $(\hat{\mathscr P},\hat{\mathscr G})[\texttt{Top}]$.
\end{definition}

The functor $R''$ induces a natural inclusion of $(\hat{\mathscr P},\hat{\mathscr G})[\texttt{Top}]$ into $\hat{\mathscr G}\wr\hat{\mathscr P}[\texttt{Top}]$. Define the functor
\begin{gather*}
    \Lambda'':\hat{\mathscr L}\wr\hat{\mathscr H}_\infty[\texttt{Top}]
    \rightarrow(\hat{\mathscr H}_\infty,\hat{\mathscr L})[\texttt{Top}],
    \quad
    \Lambda'' Z:=B(L''(\hat L\wr\hat H_\infty)R''L'',\hat L\wr\hat H_\infty,Z\!\!\restriction_{\Upsilon\wr\Upsilon}).
\end{gather*}
Elements of $\Lambda'' Z(\underline a^0_*)$ are of the form
\begin{gather*}
    [((\iota f^0,\hat\alpha^{0}),\langle (\hat f^{di_1},\langle \hat \alpha^{di_1i_{di_1}}\rangle)\rangle,\langle z^{i_1}\rangle),\langle t^d\rangle]
  =
  \\
  \left[\xymatrix@=0.75cm{
  (\underline 1_*, \underline a^{0}_* )&
    (\underline m^1_*,\langle\underline a^{1i_1}_*\rangle)
    \ar[l]_(0.5){(\iota f^0,\hat\alpha^{0})}},
    \raisebox{0.75cm}{\xymatrix@R=0.35cm@C=1.75cm{
    (\{1\}_*,\langle \underline a^{11}_*\rangle)
    &(\underline m^{(s+1)1}_*,\langle \underline a^{(s+1)1}_*\rangle)
    \ar[l]|(0.55){\ \cdots\ }_(0.55){\langle(\hat f^{d1},\langle\hat\alpha^{d1i_{d1}}\rangle)\rangle}\\
    \ar@{}@<0.3cm>[r]|(0.45){\vdots}&\\
    (\{m^1\}_*,\langle \underline a^{1m^1}_*\rangle)
    &(\underline m^{(s+1)m^1}_*,\langle \underline a^{(s+1)m^1}_*\rangle)
    \ar[l]|(0.55){\ \cdots\ }^(0.55){\langle t^d\rangle}_(0.55){\langle(\hat f^{dm^1},\langle\hat\alpha^{dm^1i_{dm^1}}\rangle)\rangle}
    }},
    \langle z^{i_1}\rangle\right],
\end{gather*}
and the $(\hat{\mathscr H}_\infty,\hat{\mathscr L})$-space structure is given by the formulas
\begin{gather*}
    \hat\alpha \langle[((\iota f^0,\hat\beta^{0}),\langle (\hat f^{di_1},\langle \hat \beta^{di_1i_{di_1}}\rangle)\rangle,\langle z^{i_1}\rangle),\langle t^d\rangle]\rangle
    \!:=\!
    [((\iota f^0,\hat\alpha\hat\beta^{0}),\langle (\hat f^{di_1},\langle \hat \beta^{di_1i_{di_1}}\rangle)\rangle,\langle z^{i_1}\rangle),\langle t^d\rangle],\\
    (f,\kappa) 
    \langle[((\iota g^{i0},\hat\alpha^{i0}),\langle (\hat g^{idj_{i1}},\langle \hat \alpha^{idj_{i1}j_{idj_{i1}}}\rangle)\rangle,\langle z^{i j_{i1}}\rangle),\langle t^{id}\rangle]\rangle
    :=\\
    \left[\left(
    \!(\iota(f\langle g^{i0}\rangle),
    \!(\eta\kappa)((\iota f)\!\ltimes\!\langle\hat\alpha^{i0}\rangle)),
    \!\left\langle\! \begin{cases}
        \!(\hat g^{idj_{i1}},\langle \hat \alpha^{idj_{i1}j_{idj_{i1}}}\rangle),
        &\!\!\!\!i=i'\\
        \!\text{id}_{(\underline n^{i(\delta^{i}(i',d)+1)}_*,\langle\underline y^{i(\delta^{i}(i',d)+1)j_{i(\delta^{i}(i',d)+1)}}_*\rangle)}\!,
        &\!\!\!\!i\neq i'
    \end{cases}\right\rangle\!,
    \!\langle z^{ij_{i1}}\rangle\right)\!,
    \!\lhd_{\underline m}\langle t^{i'd}\rangle\right].
\end{gather*}

\subsection{The {$(\mathscr H_\infty,\mathscr L)$}-spaces {$\Lambda\lVert \mathfrak A\rVert$}}

We have the adjunction
\begin{gather*}
    \xymatrix{\texttt{Top}
    \ar@<0.15cm>[r]^{R'}\ar@{}[r]|\top
    &\texttt{Top}^{\Upsilon}
    \ar@<0.15cm>[l]^{L'}};\qquad
    \begin{cases}
        R'Z\underline a_*:= Z^{\underline a}\\
        R'Z\psi:=\psi\cdot
    \end{cases},\quad
    \quad L' Z':= Z'\underline 1_*.
\end{gather*}

For an operad pair $(\mathscr P,\mathscr G,\ltimes)$ the functor $R'$ induces a natural inclusion of $(\mathscr P,\mathscr G)[\texttt{Top}]$ into $(\hat{\mathscr P},\hat{\mathscr G})[\texttt{Top}]$. We define the functor
\begin{gather*}
    \Lambda':(\hat{\mathscr H}_\infty,\hat{\mathscr L})[\texttt{Top}]\rightarrow (\mathscr H_\infty,\mathscr L)[\texttt{Top}],
    \qquad
    \Lambda' Z:=B(H_\infty L',\hat H_\infty,Z).
\end{gather*}

An element of $\Lambda'Z$ is of the form
$$
    [(\alpha^0,\langle \hat \alpha^{x_1d}\rangle,\langle z^{x_1}\rangle),\langle t^d\rangle]=\left[\alpha^0,
    \raisebox{0.75cm}{\xymatrix@R=0.35cm{
    \{1\}_*
    &\underline a_*^{1(s+1)}
    \ar[l]|{\ \cdots\ }_(0.5){\langle\hat\alpha^{1d}\rangle}\\
    \ar@{}@<0.2cm>[r]|\vdots&\\
    \{x_1\}_*
    &\underline a_*^{x_1(s+1)}
    \ar[l]|{\ \cdots\ }^(0.5){\langle t^d\rangle}_(0.5){\langle\hat\alpha^{x_1d}\rangle}
}},\langle z^{x_1}\rangle\right],$$
and the $(\mathscr H_\infty,\mathscr L)$-space structure is given by the formulas
\begin{align*}
    \alpha\langle [(\beta^{x0},
    \langle \hat \beta^{xy_{x1}d}\rangle,
    \langle z^{xy_{x1}}\rangle),\langle t^{xd}\rangle] \rangle 
    &:=
    \left[\left(\alpha\langle \beta^{a0}\rangle,
    \left\langle 
    \begin{cases}
        \hat \beta^{xy_{x1}d},
        &\!\!\! x=x'\\
        \text{id}_{\underline b^{x(\delta^{x}(x',d)+1)}},
        &\!\!\! x\neq x'
    \end{cases}
    \right\rangle,
    \langle z^{xy_{x1}}\rangle\right),\lhd_{\underline a}\langle t^{x'd}\rangle\right],\\
    f\langle [(\alpha^{i0},
    \langle \hat \alpha^{ix_{i1}d}\rangle,
    \langle z^{ix_{i1}}\rangle),\langle t^{id}\rangle] \rangle 
    &:=
    \left[\left(f\ltimes\langle \alpha^{i0}\rangle,
    \left\langle (\iota f)\ltimes \left\langle
    \begin{cases}
        \hat \alpha^{ix^i_{i1}d},
        &\!\!\! i=i'\\
        \text{id}_{\underline b^{ix^i_{i1}(\delta^i(i',d)+1)}},
        &\!\!\! i\neq i'
    \end{cases}\right\rangle
    \right\rangle,
    (f,\text{id}_{\Pi_{\underline m}\underline a^{ix^i_{i1}(s^i+1)}})\langle z^{ix^i_{i1}}\rangle\right),\lhd_{\underline m}\langle t^{i'd}\rangle\right].
\end{align*}

The functor
$$\Lambda:=\Lambda'\Lambda'':\hat{\mathscr L}\wr\hat{\mathscr H}_\infty[\texttt{Top}]\rightarrow (\mathscr H_\infty,\mathscr L)[\texttt{Top}]
$$
induces an equivalence of homotopy categories.

\section{The commutative ring spectra $\widetilde K\mathfrak A$}

The recognition principle of \cite{may2006geometry,MaEinftyRingSpcsSpectra,Ma09a} shows that the $\infty$-loop space functor $\Omega^\infty$ and the $\infty$-delooping functor $B^\infty$ induce an equivalence of the homotopy categories of ringlike of $E_\infty$-ring spaces and connective commutative ring spectra. As explained in \cite{Vi20,Vi21} this is not induced by a Quillen adjunction, but by the weaker structure of a weak Quillen quasiadjunction. The proof of this theorem in particular shows that $\pi_0\Omega^\infty B^\infty Z$ is the enveloping ring of the semiring $\pi_0 Z$, while $\pi_n \Omega^\infty B^\infty Z\cong \pi_n Z$ for all $n>0$. This is a version of Grothendieck's enveloping ring construction that is compatible with higher homotopical structure.

In this section we give a brief overview of commutative ring spectra and commutative algebra spectra, and show how to construct the $\infty$-delooping functor. In \cite{Vi21} it was convenient to construct the spectra $B^\infty Z$ as spaces of equivalence classes of filtered rooted trees with decorations, but here we use a construction that is closer to the language of categories of operators.

We construct commutative algebra spectra $\widetilde k\mathfrak A$ over the commutative ring spectra $k\mathcal U_{\mathds F}$ that represent topological $K$-theory. We will show by a general homological argument that the stable homotopy groups of $\widetilde k\mathfrak A$ coincide with the nonnegative $K$-groups of $\mathfrak A$. In order to get the correct negative stable homotopy groups implied by Bott periodicity we construct algebra spectra  $\widetilde K\mathfrak A$ by localizing at appropriate Bott elements.

\subsection{The connective commutative ring spectrum $\widetilde k\mathfrak A$}

\begin{definition}
Let $\mathds U\in\texttt{Inn}$. The category $\texttt{Sp}_{\mathds U}$ of \textit{coordinate free spectra over $\mathds U$} is an $\mathscr A_{\mathds U}$-indexed sequence of pointed spaces 
$$
    Z=\langle Z_U\rangle\in \Pi_{\mathscr A_{\mathds U}}\texttt{Top}_*
$$
equipped with structural maps 
$$
    \langle\sigma_{U,V}\rangle\in \coprod_{\Sigma_{\mathscr A_{\mathds U}}\mathscr A_U}\texttt{Top}_*(Z_U,Z_V^{\mathds S^{V-U}})
$$
satisfying
$$
    \sigma_{U,U}(z,\vec 0)=z,
    \qquad 
    \sigma_{V,W}(\sigma_{U,V}(z,\vec v),\vec w)=\sigma_{U,W}(z,\vec v+\vec w).
$$
The morphism spaces are defined as
$$
    \texttt{Sp}_{\mathds U}(Y,Z):=\{\varphi=\langle\varphi_U\rangle\in\Pi_{\mathscr A}\texttt{Top}_\ast(Y_U,Z_U)\mid 
    \sigma_{U,V}(\varphi_U y,\vec v)=\varphi_V\sigma_{U,V}(y,\vec v)\}.
$$

An \textit{$\Omega$-spectrum} is a spectrum with all structural maps weak homotopy equivalences
\end{definition}

By Brown representability $\Omega$-spectra represent (co)homology theories \cite{Brown62}. If $\mathds U=\mathds R^\infty$ then we simply write $\texttt{Sp}:=\texttt{Sp}_{\mathds R^\infty}$.  For $n\in\mathds Z$ the $n$-sphere spectrum is
$$
    \mathds S^n:=\begin{cases}
        \begin{array}{lll}
            \langle(U-\mathds R^{\lvert n\rvert})_+\rangle, & \sigma_{U,V}(\vec u,\vec v):=\vec u+\text{pr}_{V-\mathds R^{\lvert n\rvert}}\vec v, & n< 0 \\
            \langle U_+\rangle, & \sigma_{U,V}(\vec u,\vec v):=\vec u+\vec v, & n=0\\
            \langle (U\oplus \mathds R^n)_+\rangle, & \sigma_{U,V}((\vec u,\vec w),\vec v):=(\vec u+\vec v,\vec w), & n> 0
        \end{array}
    \end{cases}.
$$
We use the notation $\mathds S:=\mathds S^0$ for the sphere spectrum. For $n\in\mathds Z$ the stable $n$-homotopy group of a spectrum $Z\in\texttt{Sp}$ is
$$
\pi^S_n Z:=\pi_0\texttt{Sp}(\mathds S^n,Z).
$$
We say $Z$ is \textit{connective} if $\pi^S_n Z$ is trivial for $n<0$. If $Z$ is an $\Omega$-spectrum and $n\geq 0$ then $\pi_n^S Z\cong \pi_n Z_0$.

\begin{definition}
    For each $\underline m\in\mathbb F$ the \textit{external smash $\underline m$-product} is the functor
\begin{gather*}
    \overline \wedge_{\underline m}:\Pi_{\underline m}\texttt{Sp}\rightarrow \texttt{Sp}_{\oplus_{\underline m}\mathds R^\infty};
    \quad
    \overline \wedge_{\underline m}\langle Z^i\rangle:=\left\langle\wedge_{\underline m}Z^i_{U^i}\right\rangle,
    \quad
    \sigma_{\langle U^i\rangle,\langle V^i\rangle}([z^i],\langle \vec v^i\rangle):=[\sigma_{U^i,V^i}(z^i,\vec v^i)].
\end{gather*}
\end{definition}

For $K\subset_{\text{cpct}}\mathscr L\underline m$ and $(\langle U^i\rangle,V)\in \Sigma_{\mathscr A_{\oplus_{\underline m}\mathds R^\infty}}\mathscr A_{\sum_K f\langle U^i\rangle}$ the associated \textit{Thom complex} is
$$
    TK_V^{\langle U^i\rangle}:=\Sigma_K (V- f\langle U^i\rangle)_+/_{(f,\infty)\sim (g,\infty)}\in\texttt{Top}_\ast,
$$
where $\Sigma_K (V- f\langle U^i\rangle)_+\subset K\times V_+$ is equipped with the subspaces topology. We will use the notation $\prescript{\vec v}{f}{}:=[f,\vec v]\in TK^{\langle U^x\rangle}_{V}$, so that the base point is denoted as $\prescript{\infty}{f}{}$ for any $f\in K$. Set
\begin{gather*}
    \nu\in\texttt{POSet}(\mathscr A,\mathscr A_{\oplus_A\mathds R^\infty}),
    \qquad
    \nu  U
    := \cap_K \langle f^*_a U\rangle.
\end{gather*}
\begin{definition}
    The \textit{twisted half-smash product} is the functor
\begin{gather*}
    \mathscr L \underline m \ltimes-:\texttt{Sp}_{\oplus_{\underline m}\mathds R^\infty}\rightarrow \texttt{Sp};
    \qquad
    \mathscr L \underline m \ltimes Z:=
    \left\langle\underset{K\subset_{cpct}\mathscr L \underline m }{\text{Colim }} TK^{\nu U}_U\wedge Z_{\nu U}\right\rangle,\\
    \sigma_{U,V}\left(\left[\prescript{\vec u}{f}{},z \right],\vec v\right):= \left[\prescript{\text{pr}_{V-f\nu  V}(\vec u+\vec v)}{f}{},\sigma_{\nu U,\nu V} \left[z,f\!\!\restriction_{\nu  V}^*(\vec u+\vec v)\right]\right].
\end{gather*}
\end{definition}

The monoid $\mathscr L\underline 1$ induces a monad $(\mathds L; \eta,\mu)$ on $\texttt{Sp}$ with
$$
\mathds L Z:=\mathscr L\underline 1\ltimes Z; 
\qquad \eta_U z:\left[\prescript{\vec 0}{\text{id}}{},z\right], 
\quad \mu_U\left[\prescript{\vec u}{f}{},\left[\prescript{\vec v}{g}{},z\right]\right]:=\left[\prescript{\vec u+f\vec v}{fg}{},z\right].
$$

For $(Z,\mathfrak z)\in \mathds L[\texttt{Sp}]$ we will use the notation $\prescript{\vec u}{f}{z}:=\mathfrak z_U[\prescript{\vec u}{f}{},z]$. 

\begin{definition}
    For each $\underline m\in \mathbb F$ the \textit{smash $\underline m$-product} in $\mathds L[\texttt{Sp}]$ is defined as
$$
\wedge_{\underline m}\langle Z^i\rangle:=\left\langle\mathscr L\underline m\ltimes \overline \wedge_{\underline m} Z^i/_{\left[\prescript{\vec u}{f}{},\left[\prescript{\vec v^i}{g^i}{z^i}\right]\right]\sim\left[\prescript{\vec u+f_i\vec v^i}{f\langle g^i\rangle}{},[z^i]\right]}
\right\rangle
$$
\end{definition}

We use the notation $\otimes^{\vec u}_f [z^i]:=\left[\prescript{\vec u}{f}{},[z^i]\right]\in \wedge_{\underline m}Z^i$.

The wedge product $\wedge=\wedge_{\underline 2}$ is associative and symmetric up to natural isomorphisms. The sphere spectrum is almost a unit for the smash product, in the sense that we have a natural weak equivalence
$$
\rho_Z:Z\wedge \mathds S\rightarrow Z, \qquad \rho_{Z,U}\otimes^{\vec u}_f[z,\vec v]:=\sigma_{f_1\nu U^1,U}[\prescript{\vec u}{f_1}{z},f_2\vec v].
$$

\begin{definition}
    The category of \textit{sphere modules} $\texttt{Mod}_{\mathds S}$ is the full subcategory of $\mathds L[\texttt{Sp}]$ composed of the $\mathds L$-algebras $Z$ with $\rho_Z$ an isomorphism.
\end{definition}

The functor
$$
\Sigma^{\mathds S}:=-\wedge \mathds S:\mathds L[\texttt{Sp}]\rightarrow \texttt{Mod}_{\mathds S}
$$
is a left adjoint of the inclusion of $\texttt{Mod}_{\mathds S}$ in $\mathds L[\texttt{Sp}]$. The category $\texttt{Mod}_{\mathds S}$, equipped with the smash product as tensor product and the sphere spectrum as unit, is a symmetric monoidal category. 

\begin{definition}
    The category of \textit{commutative ring spectra} $\texttt{CRingSp}$ is the category of commutative monoids in $\texttt{Mod}_{\mathds S}$.
\end{definition}

Alternatively we have a monad $(\mathds P;\eta,\mu)$ on $\mathds L[\texttt{Sp}]$ with
\begin{gather*}
    \mathds PZ:=\int^{\mathbb F^{\text{bij}}}\!\!\!\!\wedge_{\underline m}\langle Z\rangle;
    \qquad
    \eta_U z:=\left[\otimes^{\vec 0}_{\text{id}}z\right],
    \quad
    \mu_U\left[\otimes^{\vec u}_{f}\left[\otimes^{\vec v^i}_{g^i}[z^{ij}]\right]\right]:=\left[\otimes^{\vec u+ f_i\vec v^i}_{f\langle g^i\rangle}[z^{ij}]\right].
\end{gather*}

An \textit{$E_\infty$-ring spectrum} is a $\mathds P$-algebra in $\mathds L[\texttt{Sp}]$. The monad $\mathds P$ induces a monad on $\texttt{Mod}_{\mathds S}$, such that we have an isomorphism $\texttt{CRingSp}\cong \mathds P[\texttt{Mod}_{\mathds S}]$. The functor $\Sigma^{\mathds S}$ induces a functor from $\mathds P[\mathds L[\texttt{Sp}]]$ to $\texttt{CRingSp}$, which is part of a Quillen equivalence. We will use the notation $\textstyle \prod_f^{\vec u} r^i:=\mu\left[\prescript{\vec u}{f}{},[r^i]\right]$ for spectral multiplicative structures.

\begin{definition}
    For $R\in \texttt{CRingSp}$ an \textit{$R$-module} $M$ is a sphere module equipped with an associative and unital action map $\rho:M\wedge R\rightarrow M$. 
\end{definition}    
    
The category $\texttt{Mod}_R$ of $R$-modules admits a symmetric monoidal structure induced by the wedge product of sphere modules.

\begin{definition}
    The category $\texttt{CAlgSp}_{R,\text{nu}}$ of \textit{non-unitary commutative algebra spectra} is the category on not-necessarily unital commutative monoids in $\texttt{Mod}_R$.
    
    The subcategory $\texttt{CAlgSp}_{R}$ of \textit{commutative $R$-algebra spectra} is the subcategory of unitary commutative monoids in $\texttt{Mod}_R$. 
\end{definition}

Equivalently a commutative $R$-algebra is a commutative ring spectrum $A$ equipped with a ring spectra map $\eta:R\rightarrow A$. 
A similar monad $\mathds P^{\shortrightarrow}$ on pairs of $\mathds L$-spectra can also be used to describe algebra spectra \cite{Vi21}. 

For all $R\in \mathds P[\mathds L[\texttt{Sp}]]$ the stable homotopy groups $\pi^S_- R$ form a graded commutative ring. If $R$ is connective then multiplication in $\pi^S_- R$ is given by the formula
$$\textstyle
        \prod_{\underline m}[\varphi^i]:=\left[(\vec u,\vec v)\mapsto \prod_f^{\text{pr}_{f^\perp}\vec u}\varphi^i(f^*_i\vec u, \iota_i^*\vec v)\right]
$$
for any choice of $f\in\mathscr L\underline m$, where $\langle \iota_i\rangle:\oplus_{\underline m}\mathds R^{\deg \varphi^i}\rightarrow \mathds R^{\sum_{\underline m}\deg \varphi^i}$ denotes the obvious isomorphism. Since $\mathscr L$ is an $E_\infty$-operad, i.e. equivariantly contractible, the multiplication is well defined.

For each $Z\in (\mathscr H_\infty,\mathscr L)[\texttt{Top}]$ and $U\in \mathscr A$ the inclusion of $\mathscr H_U$ in $\mathscr H_\infty$ induces by pullback a $\mathscr H_U$-space structure on $Z$.  We define the category $\Check{\mathscr H}_U$ with the same objects as $\mathbb F$ and
\begin{gather*}
    \Check{\mathscr H}_\infty(\underline a^0,\underline a^1)
    :=\left\{\Check\alpha=(\epsilon\Check\alpha,\langle \alpha^{x_1}\rangle)\in \coprod_{\mathbb F^{\text{inj}}(\underline a^0,\underline a^1)}\prod_{\underline a^1}\mathscr C(\sigma^{-1}x_1)
    \right\};\\
    \Check\alpha^2\Check\alpha^1
    :=(\epsilon\Check \alpha^2\epsilon\Check\alpha^1,\langle\langle \alpha^{2x_2}_{\epsilon\Check\alpha^1x_0}\alpha^{1(\epsilon\Check\alpha^1 x_0)}_{x_0}\rangle\rangle),
    \qquad
    \text{id}_{\underline a}:=(\text{id}_{\underline a},\langle\textbf{id}\rangle).
\end{gather*}

For each $f\in\mathscr L\underline m$ we have an action $f\ltimes: \check{\mathscr H}_U^{\underline m}\rightarrow \check{\mathscr H}_U$ defined as
$$
f\ltimes\langle \check \alpha^{i}\rangle=\langle (\wedge_{\underline m}\epsilon \check\alpha^{i},\langle f\ltimes \langle \alpha^{i x_{i1}}\rangle\rangle) \rangle.
$$

Define the monad $\check H_\infty^+$ on $\texttt{Top}_*^{\mathbb F^{\text{inj}}}$ as
\begin{gather*}
(\check H_\infty^+ Z)\underline a^0:=\int^{\mathbb F^{\text{inj}}}\hspace{-0.35cm}\check{\mathscr H}_\infty(\underline a^1,\underline a^0)\times Z\underline a^1;\quad
\eta_{\underline a_0}z:=[\text{id}_{\underline a^0}, z], 
\quad
\mu_{\underline a^0}[\check{\alpha}^0,[\check{\alpha}^1,z]]:=[\check{\alpha}^0\check{\alpha}^1,z].
\end{gather*}

We have the adjunction
\begin{gather*}
    \xymatrix{\texttt{Top}
    \ar@<0.15cm>[r]^{R}\ar@{}[r]|\top
    &\texttt{Top}^{\mathbb F^{\text{inj}}}
    \ar@<0.15cm>[l]^{L}};\qquad
    \begin{cases}
        RZ\underline a:= Z^{\underline a}\\
        RZ\psi:=\psi\cdot
    \end{cases},\quad
    \quad L Z':= Z'\underline 1.
\end{gather*}

\begin{definition}
    The \textit{$\infty$-delooping functor} is
\begin{gather*}
    B^\infty:(\mathscr H_\infty,\mathscr L)[\texttt{Top}]\rightarrow \mathds P[\mathds L[\texttt{Sp}]]; 
    \qquad 
    B^\infty Z:= B(\Sigma^UL,\check{H}^+_U,RZ).
\end{gather*}

An element of $B^\infty Z$ is of the form
$$
    [(\vec u,\langle \Check\alpha^{d}\rangle,\langle z^{x_{s+1}} \rangle),\langle t^d\rangle]=\left[\xymatrix{
    \vec u,\underline 1&
    \underline a^{s+1}
    \ar[l]^(0.475){\langle t^d\rangle}_(0.475){\langle\Check\alpha^d\rangle}|(0.475){\ \cdots\ }
},\langle z^{x_{s+1}} \rangle,\right].
$$
The inclusions $\iota_{U,V}:\mathscr H_U\hookrightarrow \mathscr H_V$ induce inclusions $\iota_{U,V}:\check{\mathscr H}_U\hookrightarrow \check{\mathscr H}_V$, and the spectra structural maps and spectral multiplication are
\begin{align*}
\sigma_{U,V}([(\vec u,\langle \Check\alpha^{d}\rangle,\langle z^{x_{s+1}} \rangle),\langle t^d\rangle],\vec v)
&:=[(\vec u+\vec v,\langle \iota_{U,V}\Check\alpha^{d}\rangle,\langle z^{x_{s+1}} \rangle),\langle t^d\rangle];\\
\textstyle
\prod^{\vec u}_f\langle[(\vec v^i,\langle \Check\alpha^{id}\rangle,\langle z^{ix_{i(s^i+1)}} \rangle),\langle t^{id}\rangle]\rangle
&:=
\left[\left(
\vec u+f_i\vec v^i,
\left\langle f\ltimes\left\langle\begin{cases}
    \check\alpha^{id},&i=i'\\
    \text{id}_{\underline a}^{i(\delta^i(i',d)+1)},&i\neq i'
\end{cases}\right\rangle\right\rangle
,f_i\langle z^{ix^i_{i(s^i+1)}}\rangle
\right),
\lhd_{\underline m}\langle t^{i'd}\rangle
\right].
\end{align*}
\end{definition}

The connective ring spectrum
$$
k\mathcal U_{\mathds F}:=
\Sigma^{\mathds S}B^\infty\Lambda\lVert\mathcal U_{\mathds F}\rVert\in \texttt{CRingSp}
$$
represents connective topological $K$-theory, real when $\mathds F=\mathds R$ and complex when $\mathds F=\mathds C$. 

Let $0_{\mathds F}$ represent the trivial $C^*$-algebra over $\mathds F$, so that $\widetilde{\mathfrak K 0_{\mathds F}}\cong\widetilde{0_{\mathds F}}\cong \mathds F$. Let
$$
     k0_{\mathds F}:=\Sigma^{\mathds S}B^\infty\Lambda\lVert 0_{\mathds F}\rVert\in \texttt{CRingSp}
$$
the natural inclusion $\mathcal U_{\mathds F}\rightarrow \texttt{pr}_{\mathds F}$ induces a homotopy equivalence
$$
k\mathcal U_{\mathds F}\xhookrightarrow{\simeq} k0_{\mathds F}.
$$

For each $\mathfrak A\in \texttt{C$^*$Alg}_{\mathds F}$ the inclusion $\lambda :\mathds F\rightarrow \widetilde{\mathfrak{KA}}$ induces a $k\mathcal U_{\mathds F}$-algebra structure on $\Sigma^{\mathds S}B^\infty\Lambda\lVert \mathfrak A\rVert$, and the projection $\pi:\widetilde{\mathfrak{KA}}\rightarrow \mathds F$ induces a natural augmentation map 
$$
    \pi_*:\Sigma^{\mathds S}B^\infty\Lambda\lVert \mathfrak A\rVert
\rightarrow
k0_{\mathds F}.
$$

Using the homotopy fiber as a homotopical analog of kernels we get connective ring spectra that represent the operator $K$-theory of $C^*$-algebras.

\begin{definition}
    The \textit{unreduced $K$-theory connective $k\mathcal U_{\mathds F}$-algebra spectra functor} is
    $$
    k:\texttt{C$^*$Alg}_{\mathds F}\rightarrow \texttt{CAlgSp}_{ k\mathcal U_{\mathds F}}, \qquad k\mathfrak A:=\Sigma^{\mathds S}B^\infty\Lambda\lVert \mathfrak A\rVert.
    $$
    
    The \textit{reduced $K$-theory connective $k\mathcal U_{\mathds F}$-algebra spectra functor} is the natural homotopy fiber
    $$
    \widetilde k:\texttt{C$^*$Alg}_{\mathds F}\rightarrow \texttt{CAlgSp}_{ k\mathcal U_{\mathds F},\text{nu}}, \qquad \widetilde k\mathfrak A:=k\mathfrak A\times_{ k0_{\mathds F}} k0_{\mathds F}^I.
    $$
\end{definition}

Elements of $k\mathfrak A$ are of the form
\begin{gather*}
    \otimes^{\vec u}_f\left[\left[\vec v^1, \mqty{\langle \Check\alpha^{d}\rangle \\ \langle t^d\rangle} ,\beta^{x_{s+1}0}
    ,\mqty{\langle \hat\beta^{x_{s+1}x'_1d'}\rangle\\\langle t^{x_{s+1}d'}\rangle},
    (\iota f^{x_{s+1}x'_10},\hat\gamma^{x_{s+1}x'_10}),
    \mqty{\langle(\hat f^{x_{s+1}x'_1d''i''_1},\langle\hat\gamma^{x_{s+1}x'_1d''i''_1i''_{d''i''_1}}\rangle)\rangle\\\langle t^{x_{s+1}x'_1d''}\rangle},\right.\right.\\
    \left.\left.
    (\hat g^{x_{s+1}x'_1i''_10},\langle\hat\kappa^{x_{s+1}x'_1i''_10i'''_0}\rangle),
    \mqty{\langle(\hat g^{x_{s+1}x'_1i''_1d'''},\langle\kappa^{x_{s+1}x'_1i''_1d'''i'''_{d'''}}\rangle,\boldsymbol U^{x_{s+1}x'_1i''_1d'''i'''_{d'''}x'''_{d'''}})\rangle\\\langle t^{d'''}\rangle},\langle \boldsymbol P^{x_{s+1}x'_1i''_1i'''_{s'''+1}x_{s'''+1}}\rangle\right],\vec v^2\right],
\end{gather*}
and elements of $\widetilde k\mathfrak A$ come further equipped with a path from the image of $\pi_*$ to the base point in $k0_{\mathds F}$.

\begin{theorem}\label{ConnectiveMainResult}
    There  is a natural isomorphism between the $\mathds N$-graded homology theories $K_-$ and $\pi^S_-\widetilde k$.
    
    The spectral multiplicative structure of $\widetilde k\mathfrak A$ induces a graded-commutative ring structure on $\widetilde k_-\mathfrak A$.
\end{theorem}

\textbf{Proof:} The natural function
    \begin{gather*}
        \Phi_{\mathfrak A}:
        K_{00}\widetilde{\mathfrak A}
        \rightarrow
        \texttt{Sp}(\mathds S,k\mathfrak A), \\
        \Phi_{\mathfrak A}([P]_0-[Q]_0)_U\vec u
        :=
        \otimes^{\text{pr}_{f^\perp}\vec u}_f\left[\left[\text{pr}_{\langle\vec e_i\mid i\geq 2\rangle}f_1^*\vec u-\ln\lvert \vec e_1\cdot f_1^*\vec u\rvert\vec e_1,
        \textbf{id},
        \text{id}_{(\underline 1_*,\underline 1_*)},
        \text{id}_{(\underline 1_*,\underline 1_*)},
        \begin{cases}
        P,&\vec e_1\cdot f_1^*\vec u\leq 0\\
        Q,&\vec e_1\cdot f_1^*\vec u\geq 0
        \end{cases}\right],f_2^*\vec u\right],
    \end{gather*}
where $\langle \vec e_i\mid i\geq 1\rangle$ is the canonical base of $\mathds R^\infty$, induces a natural homomorphism $\Phi_{\mathfrak A,*}:K_{00}\widetilde{\mathfrak A} \rightarrow
\pi^S_0k\mathfrak A$. It is a natural isomorphism, since by the recognition theorem and the fact that $\Sigma^{\mathds S}$ is part of a Quillen equivalence $\pi_0^S k\mathfrak A$ is the group completion of $\pi_0\Lambda\lVert \mathfrak A\rVert$. It induces a natural isomorphism of short exact sequences
\begin{gather*}
    \xymatrix@R=0.5cm{
        K_0\mathfrak A
       \ar@{^{(}->}[r]
       \ar[d]_{\Phi_{\mathfrak A}\restriction}^{\cong}
       &K_{00}\widetilde{\mathfrak{A}}
        \ar@{->>}[r]
        \ar[d]_{\Phi_{\mathfrak A}}^\cong
        &K_{00}\mathds F
        \ar[d]_{\Phi_{0}}^\cong\\
        \pi^S_0 \widetilde k\mathfrak A
        \ar@{_{(}->}[r]
        &\pi^S_0 k\mathfrak A
        \ar@{->>}[r]&\pi^S_0  k0_{\mathds F}
    }.
\end{gather*}
Therefore the homology theories $K_-$ and $\pi^S_- k$ are derived from naturally isomorphic functors. The second statement follows from the first. $\blacksquare$\\

The elements $[P]_0-[Q]_0\in K_n\mathfrak A$ are mapped to $\left[\left\langle
(\vec u,\vec v)
\mapsto
\Phi_{\mathfrak A}([P\vec v]_0-[Q\vec v]_0)_U\vec u
\right\rangle\right]$.

\subsection{Bott periodicity and localization}

As explained in \cite{atiyah1968bott,MaEinftyRingSpcsSpectra,EKMM} the isomorphisms of topological $K$-theoretical Bott periodicity are induced by right multiplication of virtual bundles over spheres, known as Bott elements. The same is true in operator $K$-theory, and as in the topological setting we can localize the ring spectra $\widetilde k\mathfrak A$ at these elements in order to define periodic ring spectra with the correct negative stable homotopy groups. For simplicity we only give details of localizations at Bott elements in the complex case, and just give a sketch of how the argument need to be adapted in the real case. In the notation of the present article the complex Bott element is
\begin{gather*}
    [\beta_{\mathds C}]\in \pi_2^S k 0_{\mathds C},\qquad
    \beta_{\mathds C}:=\left\langle (\vec u,\vec v)
    \mapsto
    \Phi_{0_{\mathds F}}\left(\left[\widetilde 1_{\underline 1}\right]_0-
    \left[\mqty[\frac{\lvert \vec v\rvert}{\lvert \vec v\rvert^2+1}\tilde 1
        &\frac{v_1+v_2i}{\lvert \vec v\rvert^2+1}\tilde 1\\
        \frac{v_1-v_2i}{\lvert \vec v\rvert^2+1}\tilde 1&\frac{1}{\lvert \vec v\rvert^2+1}\tilde 1]\right]_0\right)_U\vec u
    \right\rangle.
\end{gather*}

This element is associated via Swan's theorem to the virtual bundle over the 2-sphere given by the formal difference between the trivial line bundle and the tautological line bundle. The real Bott element $[\beta_{\mathds R}]\in \pi^S_8 k0_{\mathds R}$ is similarly constructed, and it is associated with the virtual line bundle over the 8-sphere given by the formal difference of a summand of the spin bundle and the trivial line bundle.

Since all the $\widetilde k\mathfrak A$ are $k0_{\mathds F}$-modules we have a well defined right multiplication map
\begin{gather*}
    R_{\beta_\mathds C}\in\texttt{Mod}_{\mathds S}(\widetilde k\mathfrak A\wedge \mathds L\mathds S^2,\widetilde k\mathfrak A)\\
    R_{\beta_{\mathds C}}\left(\otimes^{\vec u}_h\left[\left[\vec v^1, \mqty{\langle \Check\alpha^{d}\rangle \\ \langle t^d\rangle} ,\beta^{x_{s+1}0}
    ,\mqty{\langle \hat\beta^{x_{s+1}x'_1d'}\rangle\\\langle t^{x_{s+1}d'}\rangle},
    (\iota f^{x_{s+1}x'_10},\hat\gamma^{x_{s+1}x'_10}),
    \mqty{\langle(\hat f^{x_{s+1}x'_1d''i''_1},\langle\hat\gamma^{x_{s+1}x'_1d''i''_1i''_{d''i''_1}}\rangle)\rangle\\\langle t^{x_{s+1}x'_1d''}\rangle},
    \right.\right.\right.\\\left.\left.\left.
    (\hat g^{x_{s+1}x'_1i''_10},\langle\hat\kappa^{x_{s+1}x'_1i''_10i'''_0}\rangle),
    \mqty{\langle(\hat g^{x_{s+1}x'_1i''_1d'''},\langle\kappa^{x_{s+1}x'_1i''_1d'''i'''_{d'''}}\rangle,\boldsymbol U^{x_{s+1}x'_1i''_1d'''i'''_{d'''}x'''_{d'''}})\rangle\\\langle t^{d'''}\rangle},\langle \boldsymbol P^{x_{s+1}x'_1i''_1i'''_{s'''+1}x_{s'''+1}}\rangle\right],\vec v^2,[\prescript{\vec x^3}{k}{},(\vec v^3,\vec w^3)]\right]\right):=\\
    \otimes^{\vec u+h_3(\vec x^3+k\vec v^3)}_{\langle h_1,h_2\rangle}\left[\left[\text{pr}_{\langle\vec e_i\mid i\geq 2\rangle}\vec v^1-\ln\lvert \vec e_1\cdot \vec v^1\rvert\vec e_1,
    \mqty{\langle \Check\alpha^{d}\rangle \\ \langle t^d\rangle} ,\beta^{x_{s+1}0}
    ,\mqty{\langle \hat\beta^{x_{s+1}x'_1d'}\rangle\\\langle t^{x_{s+1}d'}\rangle},
    (\iota f^{x_{s+1}x'_10},\hat\gamma^{x_{s+1}x'_10}),
    \mqty{\langle(\hat f^{x_{s+1}x'_1d''i''_1},\langle\hat\gamma^{x_{s+1}x'_1d''i''_1i''_{d''i''_1}}\rangle)\rangle\\\langle t^{x_{s+1}x'_1d''}\rangle},
    \right.\right.\\
    (\hat g^{x_{s+1}x'_1i''_10},\langle\hat\kappa^{x_{s+1}x'_1i''_10i'''_0}\rangle),
    \mqty{\langle(\hat g^{x_{s+1}x'_1i''_1d'''},\langle\kappa^{x_{s+1}x'_1i''_1d'''i'''_{d'''}}\rangle,\boldsymbol U^{x_{s+1}x'_1i''_1d'''i'''_{d'''}x'''_{d'''}})\rangle\\\langle t^{d'''}\rangle},\\
    \left.\left.
    \begin{cases}
        \langle \boldsymbol P^{x_{s+1}x'_1i''_1i'''_{s'''+1}x_{s'''+1}}\rangle,&\vec e_1\cdot f_1^*\vec u\leq 0
        \\
        \left\langle \boldsymbol P^{x_{s+1}x'_1i''_1i'''_{s'''+1}x_{s'''+1}}\otimes \mqty[\frac{\lvert \vec w^3\rvert^2}{\lvert \vec w^3\rvert^2+1}\tilde 1
        &\frac{w^3_1+w^3_2i}{\lvert \vec w^3\rvert^2+1}\tilde 1\\
        \frac{w^3_1-w^3_2i}{\lvert \vec w^3\rvert^2+1}\tilde 1&\frac{1}{\lvert \vec w^3\rvert^2+1}\tilde 1]\right\rangle
        ,&\vec e_1\cdot f_1^*\vec u\geq 0
    \end{cases}\right],\vec v^2\right]
\end{gather*}

In the real case we have a similar map $R_{\beta_\mathds R}\in\texttt{Mod}_{\mathds S}(\widetilde k\mathfrak A\wedge \mathds L\mathds S^8,\widetilde k\mathfrak A)$. We have natural systems of isomorphisms $\chi_{n_1,n_2}\in \texttt{Mod}_{\mathds S}(\mathds S\wedge
\mathds L\mathds S^{\underline{n_1}}\wedge \mathds L\mathds S^{\underline{n_2}},\mathds S\wedge\mathds L\mathds S^{\underline{n_1+n_2}})$ for all $n_1,n_2\in \mathds Z$ (see \cite[III.3.7]{EKMM} and \cite[pp 386-389]{gaunce2006equivariant}). This means that we can define the map $\tilde \beta_{\mathds C}$ such that the diagram bellow commutes:
    $$
    \xymatrix@=0.5cm{
    \widetilde k\mathfrak A\wedge \mathds L\mathds S^0
    \ar[rr]^{\tilde \beta_{\mathds C}}
    \ar[dr]_(0.35){\widetilde k\mathfrak A\wedge \chi^{-1}_{2,-2}}&&\widetilde k\mathfrak A\wedge \mathds L\mathds S^{-2}\\
    &\widetilde k\mathfrak A\wedge \mathds L\mathds S^{2}\wedge \mathds L\mathds S^{-2}
    \ar[ur]_(0.65){R_{\beta_{\mathds C}}\wedge\mathds L\mathds S^{-2}}&
    }
    $$

We can repeatedly apply $\tilde \beta_{\mathds C}$ to get an infinite sequence
$$
\xymatrix@=1.25cm{
\widetilde k\mathfrak A\wedge \mathds L\mathds S^0
\ar[r]^{\tilde \beta_{\mathds C}}&
\widetilde k\mathfrak A\wedge \mathds L\mathds S^{-2}
\ar[r]^{\tilde \beta_{\mathds C}\wedge \mathds L\mathds S^{-2}}&
\widetilde k\mathfrak A\wedge \mathds L\mathds S^{-4}
\ar[r]^(0.55){\tilde \beta_{\mathds C}\wedge \mathds L\mathds S^{-4}}&\cdots
}.
$$

The localization at the Bott element is then the telescope construction on this sequence, i.e.
$$
    \widetilde k \mathfrak A[\beta_{\mathds C}^{-1}]
    :=\text{Tel}(\tilde\beta_{\mathds C}\wedge \mathds L\mathds S^{-2\bullet})
    =\left(\coprod_{\mathds N}\widetilde k\mathfrak A\wedge \mathds{LS}^{-2n}\wedge [n,n+1]_+\right) /_\sim
$$
with $\widetilde k\mathfrak A\wedge \mathds{LS}^{-2n}\wedge\{n+1\}_+$ being glued to $\widetilde k\mathfrak A\wedge \mathds{LS}^{-2n-2}\wedge\{n\}_+$ by the map $\tilde\beta_{\mathds C}$. We have a similar localization $\widetilde k \mathfrak A[\beta_{\mathds R}^{-1}]$ in the real case induced by $R_{\beta_{\mathds R}}$. 

We may also localize the ring spectra $k\mathcal U_{\mathds F}$ at Bott elements, giving us the commutative ring spectra 
$$
    K\mathcal U_{\mathds F}:=k\mathcal U_{\mathds F}[\beta_{\mathds F}^{-1}]
$$
that represents periodic topological $K$-theory. We then have that $\widetilde k\mathfrak A[\beta_{\mathds F}^{-1}]$ are $K\mathcal U_{\mathds F}$-algebras.

\begin{definition}
    The $K$-theory $k\mathcal U_{\mathds F}$-algebra spectra functor is
    $$
    \widetilde K:\texttt{C$^*$Alg}_{\mathds F}\rightarrow \texttt{CAlgSp}_{ K\mathcal U_{\mathds F},\text{nu}}, \qquad \widetilde K\mathfrak A:=\widetilde k\mathfrak A[\beta_{\mathds F}^{-1}].
    $$
\end{definition}

By construction of the localization and theorem \ref{ConnectiveMainResult} we get our main result.

\begin{theorem}
    There  is a natural isomorphism between the $\mathds Z$-graded homology theories $K_-$ and $\pi^S_-\widetilde K$.
    
    The spectral multiplicative structure of $\widetilde K\mathfrak A$ induces a graded-commutative ring structure on $K_-\mathfrak A$.
\end{theorem}

\bibliographystyle{plain}
\bibliography{bibliography.bib}

\end{document}